\documentclass[12pt]{amsbook}
\usepackage{amssymb}
\usepackage[T1]{fontenc}
\def\l{\lambda}

\def\t{\times}

\newtheorem{theo}{Theorem}[section]
\newtheorem{prop}[theo]{Proposition}
\newtheorem{lemme}[theo]{Lemma}
\newtheorem{cor}[theo]{Corollary}

\def\vol{{\rm vol}}

\def\cc{{\mathbb C}}
\def\n{{\mathbb N}}
\def\z{{\mathbb Z}}
\def\q{{\mathbb Q}}
\def\r{{\mathbb R}}
\def\a{{\mathbb A}}
\def\bc{\backslash}

\def\bc{\backslash}

\def\tr{{\rm tr}}

\def\o{{\omega}}

\def\ind{{\rm ind}}
\def\g{\gamma}

\def\h{{\mathbb H}}

\def\H{{\mathcal H}}

\def\lra{\leftrightarrow}

\begin{document}

\centerline{\Large{\bf The trace formula}}
\ \\

\centerline{\Large{\bf and the proof of the global Jacquet-Langlands correspondence}}
\ \\

\centerline{A.I.Badulescu}
\centerline{IMAG, Univ Montpellier, CNRS, Montpellier, France}

\newpage

\tableofcontents

\section{Introduction}

This paper contains the material covered in the lectures I gave at the doctoral school {\it Introduction to Relative Aspects in Representation Theory, Langlands Functoriality and Automorphic Forms} held at the CIRM, Luminy, 16-20 May 2016. The organizers, Dipendra Prasad, Volker Heiermann and Fiona Murnaghan, asked me to give lectures on the proof of the global Jacquet-Langlands correspondence, following my paper [Ba1], as an application of the trace formula of Arthur. The lectures were designed to meet the level of PhD students. I was very pleased to receive the invitation and I would like to warmly thank the organizers for this opportunity and for the perfect organization and wonderful atmosphere at the conference. 

I was a bit scared however when  they asked me to write a paper about it. In a lecture one may give an overview without getting into details, and in this field, when entering into definitions and details, everything becomes so complicated that the reader might easily give up. For example, one may say at the blackboard "by some kind of linear independence of characters we may show that...", but that looks awkward in a paper, where the reader would expect to have an explanation about which kind of linear independence of characters is actually used, since obviously it is not the classical one. Such kind of casual claim may involve new definitions and take two pages when explained - as a matter of fact, this claim became here an Appendix. 

In the end I decided to first give a fake proof, explaining what are (in my opinion, following mainly [JL] and [DKV] where I learned this material) the six standard steps of the proof of a correspondence. Fake proof based on some simplifications (for example assuming that every reductive group has a simple trace formula, which is not the case), in order for the beginner to see why such and such object appears in the field, how such and such question arose, and also because the proof of correspondences is beautiful, but its beauty is often hidden by the very complicated technical aspects. That is why I decided to postpone these aspects. The truth is that each of this steps is awfully complicated for general groups, and even for any group which is not compact (modulo its center) like $GL_n$. Subsequently, I explain step by step what was untrue in my assumption and how things go in real life. It is a way of putting all the details after, and not cut the story of the proof every five lines. In the last section, I explain the proof of the global Jacquet-Langlands correspondence following the six steps. It is a correspondence between the discrete series of the group $GL_n$ and the discrete series of an inner form of $GL_n$, over an ad\`ele ring of characteristic zero. It would have been impossible even to sketch the proof of the step 2 and part of 3, because it is extremely long and complicated. It has been carried out by Arthur and Arthur-Clozel; I will directly use their final result about the comparison of the trace formulae in [AC]. 

The paper is divided into two chapters and an appendix. In the first chapter, the reader will find a collection of results concerning restricted products, construction of ad\`ele groups of reductive groups over number fields, Hecke algebras and the theory of admissible representations of reductive groups over local fields including parabolic induction, definition of discrete series (which are representations of ad\`ele groups of reductive groups) and a proof of the trace formula in the compact quotient case. It has the ambition to give a clear overview of these matters, without proofs but with references, on a master level. I restrain myself and discuss the zero characteristic only. In the second chapter, I study general linear groups over finite dimension division algebras over local or number fields. I use the definitions from the first chapter and explain the transfer between general linear groups, then define the local Jacquet-Langlands correspondence for unitary representations (all these things are necessary in order to be able to state the global Jacquet-Langlands correspondence). Afterwards, I explain the interaction between local Jacquet-Langlands correspondence and parabolic induction, which is necessary in the proof of the global Jacquet-Langlands correspondence, where we need to use induction on the rank of the group and so proper Levi subgroups come into play. Only then I state  the global Jacquet-Langlands correspondence. Next, I give the six steps of the proof, using some simplifications and vague explanations in order to make  it less difficult to swallow. I further clarify by explaining step by step what are the simplifications I made, and give the example of $GL_2$ (the original case in [JL], here I follow [GJ]). Finally, I come to the ultimate goal of the paper and give more precise technical details of the proof of the global Jacquet-Langlands correspondence for $GL_n$, following the six steps scheme. In the appendix, I explain an important step of the proof, the classical simplification of the equality of the spectral sides (here this is "Step 4"). I give the proof for a general reductive group. 

I want to thank Abderrazak Bouaziz for useful discussions and Guy Henniart who made useful remarks on the manuscript. I also want to thank Radu Craiu for remarks on the text.

\chapter{Reductive groups}

\section{Restricted products of groups}\label{restprod}

I follow chapter 5 of [RV]. Let $V$ be a non empty set. {\it Almost all} $v\in V$ means {\it all but a finite number of elements $v$ of $V$}. If $(G_v)_{v\in V}$ is a family of groups, we denote by $\prod_{v\in V}G_v$ the product group of the $G_v$. If $g\in \prod_{v\in V}G_v$ we denote by $g_v$ the projection of $g$ on $G_v$ and call $g_v$ the local component of $g$ at $v$. Sometimes we write $(g_v)_{v\in V}$ for $g$.

Let $(G_v)_{v\in V}$ be a family of locally compact groups. Let $V_\infty $ be a finite subset of $V$, and assume for every $v\in V\bc V_\infty $ an open compact subgroup $K_v$ of $G_v$ is given. Then the {\bf restricted product} of the $G_v$, $v\in V$, with respect to the $K_v$, $v\in V\bc V_\infty$, is by definition the following topological group $G$:

- as a group, $G$ is given by
$$G:=\{g\in \prod_{v\in V} G_v,\ g_v\in K_v \ \ {\text{for almost all}}\ v\in V \}$$
with the componentwise product;

- the topology on $G$ is defined by the fact that the following sets form a neighborhood base of $1$: the products $\prod_{v\in V} N_v$, where $N_v$ is a neighborhood of $1$ in $G_v$ and, for almost all $v\in V\bc V_\infty$, $N_v=K_v$.\\
\ \\
{\bf Example.} Set $\q_\infty:=\r$. If $P$ is the set of prime numbers and $V_\infty :=\{\infty\}$, if $V:=P\cup V_\infty $, then one may define the restricted product $G$ of the $\q_p$, $p\in V$, with respect to the $\z_p$, $p\in V\bc\{\infty\}$. Then $G$ is usually denoted by $\a$ and is the {\bf ring of adèles} of $\q$ (beside being a group, it is obviously a ring). The  {\bf group of idèles} of $\q$ is the restricted product of the $\q_p^\t$, $p\in V$, with respect to the $\z_p^\t$, $p\in V\bc V_\infty$. We denote it by $\a^\t$, because as a group it is the group of invertible elements of $\a$; but the topology of $\a^\t$ is not the restriction of the topology of $\a$.

If $x\in\q$, then $(g_v)_{v\in V}$ defined by $g_v=x$ for all $v\in V$ defines an element $a_x$ of $\a$. The map $x\mapsto a_x$ is the {\bf standard embedding} of $\q$ as a subring of $\a$, and of $\q^\t$ as a subgroup of $\a^\t$.\\
\ \\
{\bf Notation.} If $S$ is a finite subset of $V$, we set $G_{S}:=\prod_{v\in S} G_v$ and we let $G^{S}$ be the restricted product of the groups $G_v$, $v\in V\bc S$, with respect to the $K_v$, $v\in V\bc (S\cup V_\infty)$.\\

For every finite subset $S$ of $V$ containing $V_\infty$ let $G(S):=\prod_{v\in S}G_v \t \prod_{v\in V\bc S}K_v$ with the product topology. Then the restricted product of the $G_v$, $v\in V$, with respect to the $K_v$, $v\in V\bc V_\infty$, is the inductive limit of the  topological subgroups $G(S)$ where $S$ runs through the finite subsets of $V$ containing $V_\infty$.
It follows that, if $S$ is a finite subset of $V$ containing $V_\infty $, then the restricted product of $G_v$, $v\in V$, with respect to the $K_v$, $v\in V\bc S$, equals $G$. Moreover, if $S$ is a finite subset of $V$, $G$ is canonically isomorphic to the direct product 
$$G_{S}\t G^{S}.$$
The isomorphism is a homeomorphism for the following (natural) topologies:  

- the topology on $G$ is the restricted product topology with respect to the $K_v$, $v\in V\bc V_\infty$,

- the topology of $G_S\t G^{S}$ is the product topology, where 

- the topology on $G_S=\prod_{v\in S}G_v$ is the (finite) product topology and 

- the topology on $G^S$ is the restricted product topology with respect to the $K_v$, $v\in V\bc S$.\\
We have the smooth embedding of $G_S$ into $G$ given by $g\mapsto (g,1)\in G_{S}\t G^{S}$, hence a smooth embedding of $H$ into $G$ for any subgroup $H$ of $G_S$.

\begin{prop}
{\rm (a)} $G$ is locally compact.

{\rm (b)} Any compact subset of $G$ is included in a product $\prod_{v\in V}H_v$ where $H_v$ is a compact subgroup of $G_v$, and, for almost all $v\in V\bc V_\infty $, $H_v=K_v$.
\end{prop}

Follows from the discussion above. See Proposition 5-1 [RV].\\

\section{Characters}\label{restprodchar}

If $H$ is a topological group, let $X(H)$ be the group of smooth characters\footnote{Here all the characters have complex values; smooth means continuous; when $H$ is locally compact and totally disconnected, a smooth character has open kernel, because open subgroups form a basis of neighborhoods of $1$ in $H$, and there is a neighborhood of $1$ in $\cc^\t$ containing no non trivial subgroup of $\cc^\t$.} of $H$. Let $V$, $V_\infty$, $G_v$, $K_v$ and $G$ be like in the previous section. If $v\in V\bc V_\infty $, if $\chi\in X(G_v)$, we say $\chi$ is {\bf unramified} if $\chi$ is trivial on $K_v$.

\begin{prop}
Assume the groups $G_v$ are abelian.  

{\rm (a)} If $\chi\in X(G)$, then $\chi$ induces by restriction a smooth character $\chi_v$ of $G_v$ for all $v\in V$. 
For almost all $v\in V\bc V_\infty $, $\chi_v$ is unramified.

{\rm (b)} If for every $v\in V$ a smooth character $\chi_v$ of $G_v$ is given, such that, for almost all $v\in V\bc V_\infty $, $\chi_v$ is unramified, then 
$$\chi(g)=\prod_{v\in V}\chi_v(g_v)$$ 
is well defined for every $g=(g_v)_{v\in V}\in G$ and is a smooth character of $G$. Moreover, the $\chi_v$ are obtained from $\chi$ as in {\rm (a)}.
\end{prop}

See Lemmas 5-2 and 5-3 [RV]. The character $\chi$ of (b) is the {\bf product character of the $\chi_v$}.\\
\ \\
{\bf Example.} Let us go back to idèles. There are canonical characters $\chi_p\in X(\q_p^\t)$, $p\in P$, 
and $\chi_\infty\in X(\q_\infty^\t)=X(\r^\t)$, given by $\chi_p(x)=|x|_p$, where $|\ |_p$ is the normalized (i.e. the value in $p$ of $|\ |_p$ is $\frac{1}{p}$) $p$-adic norm of $x$ and  $\chi_\infty(x)=|x|$. Then all the $\chi_p$, $p\in P$, are unramified and the product $\prod_{v\in V}\chi_v$, denoted  $|\ |_\a$, is a character of $\a^\t$. It is easy to check that, given the definition of the normalized $p$-adic norm, $|\ |_\a$ is trivial on the subgroup $\q^\t$ of $\a^\t$.
This character plays an important role in arithmetics. For example, it appears in the classification of the residual spectrum of $GL_n(\a)$ ([MW1]).  

We have (Proposition 6-12 [RV]) an isomorphism of topological groups 
$$\a^\times\simeq \q^\t\times\r_+^\times\times\prod{\z_p}^\t,$$
which is also a homeomorphism. If we see $\q^\t$ as a subgroup of $\a^\t$ via the standard embedding (Section \ref{restprod}), the isomorphism sends $x$ to $(x,1,1)$ for every $x\in\q^\times$. 
It follows that   
{\it every character of $\a^\times$ trivial on $\q^\times$ is the product of a character of finite order and a complex power of $|\ |_\a$} (Proposition 6-13 [RV]).\\

\section{Measures}\label{restprodmeas}

For $v\in V$, let $dg_v$ be a left Haar measure on $G_v$, such that $vol(K_v)=1$ for almost all $v\in V\bc V_\infty $. Then there exists a unique left Haar measure on $G$ such that, for every finite subset $S$ of $V$ containing $V_\infty $, the restriction of $dg$ to $G_S$ is the product measure (Proposition 5-5 [RV]). The measure $dg$ is the {\bf product measure of the $dg_v$}. Then we have:

\begin{prop}

{\rm (a)} If $f:G\to \cc$ is an integrable function, one has 
$$\int_G f=\lim_S\int_{G_S}f_S$$
where $S$ runs over the finite subsets of $V$ and $f_S$ is the restriction of $f$ to $G_S$.

{\rm (b)} If for all $v\in V$ a smooth integrable function $f_v:G_v\to \cc$ is given such that $f_v$ is the characteristic function of $K_v$ for almost all $v\in V\bc V_\infty$, then $f$ is integrable and 
$$\int_G f=\prod_{v\in V} \int_{G_v}f_v.$$
\end{prop}

See Proposition 5-6 [RV].\\
\ \\
{\bf Example.} Measures on the ring of ad\`eles and on the group of id\`eles: on $\r$ we consider the Haar measure $dg_\infty$ given by $vol([0,1])=1$ and on $\r^\t$ the measure $dg^*_\infty:=\frac{dg_\infty}{|\ |}$. On $\q_p$, $p\in P$, we fix the Haar measure $dg_p$ such that $vol(\z_p)=1$ and on $\q_p^\t$ the measure $dg^*_p$ such that $vol(\z_p^\t)=1$. One may show that $dg^*_p:=\frac{dg_p}{\ |\ |_p}$. We obtain then product Haar measures $dg$ on $\a$ and $dg^*$ on $\a^\t$, such that 
$dg^*=\frac{dg}{\ |\ |_\a}$.\\

Using the decomposition $\a^\times\simeq \q^\t\times\r_+^\times\times\prod{\z_p}^\t$, it is clear that the volume of  $\a^\t/\q^\t$ is not finite.\\

\section{Restricted tensor products of vector spaces}\label{restprodhilb}

I follow [Fl1]. Let $V$ be a nonempty set and $V_\infty$ a finite subset of $V$. Let $\{E_v\}_{v\in V}$ be a family of vector spaces over a field $F$.
For every $v\in V\bc V_\infty$ fix a vector $e_v\in E_v\backslash \{0\}$. For any finite subset $S$ of $V$ containing $V_\infty$, let $E_S=\otimes_{i\in S} E_i$. If $S'$ is a finite subset of $V$ containing $S$, there is a unique linear map $E_S\to E_{S'}$ extending   
$e\mapsto e\bigotimes(\otimes_{_{v\in S'\backslash S}} e_v)$. Let $E$ be the injective limit of this system, where $S$ runs over the finite subsets of $V$. 
We say $E$ is the {\bf restricted tensor product of spaces} $E_v$, $v\in V$, with respect to the $e_v$, $v\in V\bc V_\infty$. $E$ is spanned by vectors of type $\otimes_{_{v\in V}} x_v$, where $x_v\in E_v$ for all $v\in V$, and $x_v=e_v$ for almost all $v\in V\bc V_\infty$.\\

\section{Restricted tensor products of algebras}\label{restprodalg}

If the $E_v$ in the preceding section are $F$-algebras, if $e_v$ is an idempotent of $E_v$ for almost all $v\in V\bc V_\infty$, then the restricted tensor product of spaces $E_v$ is an algebra for the componentwise multiplication.\\

\section{Ad\`eles of reductive groups}\label{adelgroup}

The references for this section are [Fl1] and [Ti]. Let $F$ be a number field and $G$ a connected reductive algebraic group defined over $F$. Recall a place of $F$ is an equivalence class of norms on $F$. For a place $v$ of $F$ we denote by $F_v$ the completion of $F$ with respect to $v$, and by $G_v$ the group $G(F_v)$.
Let $V$ be the set of places of $F$, $V_\infty$ the (finite) set of infinite places (i.e. such that $F_v$ is Archimedean) and $V_f$ the set of finite places of $F$ (i.e. such that $F_v$ is non Archimedean). 

For each $v\in V_f$, let $O_v$ be the ring of integers of $F_v$.
Fix an algebraic groups embedding of $G(F)$ into a matrix group $GL_n(F)$. It gives rise, for every $v\in V$, to an embedding of $G_v$ into $GL_n(F_v)$. For $v\in V_f$, set  $K_v:=G_v\cap GL_n(O_v)$. Then $K_v$ is an open compact subgroup of $G_v$. Moreover, there exists a finite set $V_0$, such that $V_\infty\subset V_0\subset V$ and, for all places $v\in V\backslash V_0$, $K_v$ is a {\it maximal} open compact subgroup of $G_v$, and even hyperspecial ([Ti], 3.9; I do not give the definition of hyperspecial here).

The group $G(\a)$, is defined as the restricted product of the $G_v$, $v\in V$, with respect to the $K_v$, $v\in V\backslash V_\infty$, and it is locally compact. It is the {\bf adèle group} of $G$ over $F$.   
We set $G_\infty :=G_{V_\infty}$ and $G_f:=G^{V_\infty}$ (see Section \ref{restprod} for the notation $G_S$ and $G^S$), and identify $G(\a)$ with the direct product $G_\infty \t G_f$. 
For $v\in V_f$, we fix on $G_v$ the unique Haar measure such that  $vol(K_v)=1$. 
We fix arbitrary Haar measures on $G_v$, $v\in V_\infty$. Because the $G_v$ are reductive, a left Haar measure is also a right Haar measure ([Re], Proposition V.5.4).
On $G(\a)$ we fix the measure associated to those local choices like in Section \ref{restprodmeas}.

The group $G(F)$ has a {\bf standard embedding} in $G(\a)$ given by $x\mapsto (g_v)_{v\in V}$ where $g_v=x$ for all $v\in V$. In order to show it, one has to check that any element of $G(F)$ belongs to $K_v$ for almost all $v\in V_f$. Because $G(F)$ is a subgroup of $GL_n(F)$, and using the definition of $K_v$, it is enough to show that given a matrix $M\in GL_n(F)$, one has $M\in GL_n(O_v)$ for almost all $v\in V_f$. But any element of $F$ belongs to $O_v$ for almost all $v\in V_f$, so for almost all $v\in V_f$ all the coefficients of $M$ and of $M^{-1}$ are in $O_v$. 

In the sequel, we will treat $G(F)$ as a subgroup of $G(\a)$. It is a discrete subgroup.\\

\section{Local Hecke algebras, representations and traces}\label{heckealg}

From now on, $F$, $V$, $G$, $G_v$ etc. are like in the previous section. 
Fix  $v\in V_f$. Let 
$$\H(G_v):=\{f:G_v\to\cc, \ f \text{ is locally constant with compact support}\}.$$ 

We will sometimes denote $\H(G_v)$ simply by $\H_v$.

Endowed with the multiplication $*$ given by the convolution of functions, $\H_v$ is an algebra, without unit if $G_v$ is not compact, called the {\bf Hecke algebra} of $G_v$. The elements $e_{v,K}:=\frac{1}{{\rm vol}(K)}1_K$ ($1_K$ is the characteristic function of $K$), for open compact subgroups $K$ of $G_v$, are idempotent elements in $\H_v$. 



If $v\in V_f$, a function in $\H_v$ is called {\bf spherical} if it is left and right invariant by $K_v$. For $v$ such that $K_v$ is a hyperspecial compact subgroup (that is the case for almost all $v\in V_f$), the subalgebra $\H(G_v//K_v)$ of $\H_v$ made of spherical functions is commutative (this is explained in [Ti] 3.3.3 when $K_v$ special; but hyperspecial implies special -- [Ti] 1.9, 1.10).\\
\ \\

If $(\pi,W)$ is a representation of $G_v$, then $\pi$ is said to be {\bf smooth} if every vector in $W$ is fixed by an open subgroup of $G_v$ and {\bf admissible} if it is smooth and for every open compact subgroup $K$ of $G_v$, the space $W^K$ of vectors in $W$ fixed by $K$ is finite-dimensional. Any  smooth irreducible representation is admissible ([Re] Th.VI.2.2). If $\pi$ is admissible, if $W$ is a Hilbert space, then $\pi$ is said to be {\bf unitary} if it preserves the Hilbert product. More generally, $\pi$ is {\bf unitarisable} if there exists a Hilbert space structure on $W$ for which $\pi$ is unitary.\\

If $\pi$ is smooth and $f\in\H_v$, then we define the operator $\pi(f):W\to W$ by 
$$\pi(f)(x):=\int_{G_v} f(g)\pi(g)(x)\ dg.$$ 
When $x$ is fixed, the integral is a finite sum.
When $\pi$ is an admissible representation and $f\in\H_v$, the operator $\pi(f)$ is of finite rank and so its trace $\tr(\pi(f))$ is defined. The map  $\tr\pi:\H_v\to\cc$ defined by $f\mapsto \tr(\pi(f))$ is a $\cc$-linear form on $\H_v$ called the {\bf distribution character of $\pi$}. It determines $\pi$ up to equivalence (as the trace in the finite group case). There is actually a stronger result ([JL] Lemma 7.1, [Re] Prop.III.1.13), which will be used later:\\
\ \\
{\bf  The linear independence of the distribution characters:} {\it if $\l_1, \l_2,...,\l_k\in \cc$ and $\pi_1,\pi_2,...,\pi_k$ are smooth irreducible non isomorphic representations of $G_v$, then $\sum_{i=1}^k\lambda_i\ \tr{\pi_i} = 0$ implies all the $\l_i$ are zero}.\\ 

Notice that $\pi(g)$ for $g\in G_v$ is not in general a trace operator, so the trace is not directly definable.\\

Let $(\pi_v,W)$ be an irreducible admissible representation of $G_v$. Let $W^{K_v}$ be the space of fixed vectors under $K_v$. Then $x\mapsto \pi(f)x$, $f\in \H(G_v//K_v)$, makes $W^{K_v}$ into an irreducible $\H(G_v//K_v)$-module. When $\H(G_v//K_v)$ is commutative (we saw that occurs for almost all $v\in V_f$), $W^{K_v}$ has dimension at most one ([Car], page 151). If the dimension of $W^{K_v}$ is one, we say $\pi_v$ is {\bf unramified}.

We have the 
\begin{lemme}\label{spherical}
If $f_v$ is spherical, then for every $x\in W$, $\pi_v(f_v)(x)\in W^{K_v}$. If $\H(G_v//K_v)$ is commutative, if $\pi_v$ is not unramified, then $\pi_v(f_v)=0$, while if $\pi_v$ is unramified then $\tr\pi_v(1_{K_v})=1$.
\end{lemme}
\ \\
{\bf Proof.} Easy.\qed
\ \\

When $v\in V_\infty$ -- so $G_v$ is a real group -- a notion of {\bf admissible} representation of $G_v$ may be defined and if $\H(G_v)$ is the algebra of $C^{\infty}$-functions on $G_v$ with compact support, then for every irreducible admissible representation $\pi$ of $G_v$ and every $f\in\H(G_v)$ the operator $\pi(f)$ is defined like in the case $v\in V_f$ and it is a trace class operator. See [Kn1], Chapter VIII and Chapter X ([Del], 4, for a summary). Notice that all the  representations appearing in our paper are {\it unitarizable}. 

\section{Parabolic induction}\label{parabind}

Let $v\in V_f$. Let $P$ be a parabolic subgroup of $G_v$ and $P=LN$ a Levi decomposition of $P$. Let $\delta_P:P\to\cc$ be the modulus character ([Cas] 3.1). 
Let $\pi$ be a smooth representation of $L$ in a space $V$. Let $W$ be the space of functions $h:G_v\to \cc$ satisfying the following conditions:

- $h(mng)=\delta_P^{-1/2}(m)\pi(m)f(g)$ for all $m\in L$, $n\in N$ and $g\in G_v$,

- there is an open subgroup $K$ of $G_v$ such that $h(gk)=h(g)$ for all $g\in G_v$ and all $k\in K$.

The representation of $G_v$ in $W$ by right translations is {\bf the induced representation from $\pi$} to $G_v$ (with respect to $P$) and is denoted by $\ind_P^{G_v}\pi$. It is a smooth representation. Moreover, if $\pi$ is admissible, $\ind_P^{G_v}\pi$ is admissible ([Re] Lemma III.2.3); if $\pi$ is unitary, $\ind_P^{G_v}\pi$ is unitarizable for a standard induced scalar product ([Re] IV.2.3 since $P\bc G_v$ is compact); if $\pi$ is of finite length, $\ind_P^{G_v}\pi$ is of finite length ([Re] Lemma VI.6.2 and Proposition VI.1.2).

The theory of parabolic induction may be developed for admissible representations of real groups $G_v$, $v\in V_\infty$, and has the same properties (see the nice paper of Delorme [Del] for a summary of the theory, with precise references to [Kn1]).\\

\section{Global Hecke algebras}\label{globalheckealg}

Let $\H(G(\a))$ be the restricted tensor product of the $\H_v$, $v\in V$, with respect to the idempotents $1_{K_v}$ (the characteristic function of $K_v$), $v\in V_f$ (Section \ref{restprodalg}). 

A {\bf simple tensor} in $\H(G(\a))$ is a function $f=\otimes_{v\in V} f_v$, where $f_v\in \H_v$ for all $v\in V$ and $f_v=1_{K_v}$ for almost all $v\in V$. Any element of $\H(G(\a))$ is a linear combination of simple tensors.\\

\section{$L^2(G(F)Z(\a)\bc G(\a))$ and discrete series}\label{discreteseries}

In this section and the following, "almost all $x$ ..." means "all $x$ outside a negligible (i.e. measure zero) set". 
Let $G$ be like before and denote by $Z$ the center of $G$ and by $Z_v$ the center of $G_v$ for all $v\in V$. Then $Z(\a)$ is the restricted product of the $Z_v$, $v\in V$, with respect to the $K_v\cap Z_v$, $v\in V_f$. Fix Haar measures on $Z_v$, $v\in V$, such that, if $v\in V_f$, then the volume of $K_v\cap Z_v$ is one. Put on $Z(\a)$ the product measure $dz$ as in Section \ref{restprodmeas}. Then $Z(\a)\bc G(\a)$ inherits a quotient measure. The group $Z(\a)\bc G(\a)$ is the restricted product of the $Z_v\bc G_v$, $v\in V$, with respect to the $(K_v\cap Z_v)\bc K_v$, $v\in V_f$, and it is locally compact.

$Z(F)\bc G(F)$ is a discrete subgroup of $Z(\a)\bc G(\a)$, and we consider the standard quotient measure $dg$ on the quotient 
$G(F)Z(\a)\bc G(\a)$. The volume of $G(F)Z(\a)\bc G(\a)$ is then finite ([PR] Theorem 4.14 combined with the Theorem 5.5 and its proof; see also the original work [Bo], [Go2]). 

\begin{lemme}\label{negligible}
Let $p:G(\a)\to G(F)Z(\a)\bc G(\a)$ be the canonical projection. Let $N$ be a measurable subset of $G(F)Z(\a)\bc G(\a)$. Then $N$ is negligible if and only if $p^{-1}(N)$ is a negligible subset of $G(\a)$. 
\end{lemme}
\ \\
{\bf Proof.} The same result with "locally negligible" instead of "negligible" is the Proposition 6 (a), [Bou1] Chapter VII, section 2, subsection 3. Locally negligible sets are treated in [Bou2] Chapter IV, section 5, subsection 2, and is proved (Corollary 3) that, in locally compact spaces countable at infinity, a set is locally negligible if and only if it is negligible.\qed
\ \\

Fix a smooth unitary character $\o$ of $Z(F)\bc Z(\a)$ (see Section \ref{restprodchar}); $\o$ will be often seen as a unitary character of $Z(\a)$ trivial on $Z(F)$. 

Let $X_\o$ be the set of measurable functions $\phi :G(\a)\to\cc$ such that 

-  $\phi(zg)=\o(z)\phi(g)$ for all $z\in Z(\a)$ and all $g\in G(\a)$,  

-  $\phi(ug)=\phi(g)$ for all $u\in G(F)$ and all $g\in G(\a)$.

If $\phi,\phi'\in X_\o$, then the relation $\phi(x)=\phi'(x)$ is stable by left multiplication with elements of $G(F)Z(\a)$. We let $L^2_\o(G(F)Z(\a)\bc G(\a))$ be the space of classes of functions $\phi\in X_\o$ modulo functions which are zero for almost all $x\in G(F)Z(\a)\bc G(\a)$, such that

- $|\phi|^2$ (which is trivial on $G(F)Z(\a)$) is integrable over $G(F)Z(\a)\bc G(\a)$.\\
\ \\
The product $(\phi,\psi)=\int_{G(F)Z(\a)\bc G(\a)}\bar{\phi}\psi$ is well defined and makes $L^2_\o(G(F)Z(\a)\bc G(\a))$ into a Hilbert space.\\

From now on we fix $\o$ and we write $L^2$ instead of $L^2_\o$ for the sake of simplification.\\

Let $L^2_0(G(F)Z(\a)\bc G(\a))$ be the subspace of $L^2(G(F)Z(\a)\bc G(\a))$ made of classes of functions $\phi$ such that, for almost all $x\in G(\a)$, 
$$\int_{N(F)\backslash N(\a)} \phi(nx)dn=0$$
for every unipotent radical $N$ of a proper parabolic subgroup of $G$ defined over $F$. Notice that, if $x$ verifies this condition, then every element of $G(F)Z(\a)x$ does. Thanks to the Lemma \ref{negligible} one may say "for almost all $x\in G(F)Z(\a)\bc G(\a)$" instead of "for almost all $x$ in $G(\a)$".\\

The group $G(\a)$ acts on $L^2(G(F)Z(\a)\bc G(\a))$ by right translations. We call this representation $R$. Then $R$ is a unitary representation. 
I will call here {\bf irreducible} subrepresentation of $R$ a topological irreducible subrepresentation, i.e. a subrepresentation of $R$ which does not contain any proper {\it closed} subspace.
The representation $R$ stabilizes $L_0^2(G(F)Z(\a)\bc G(\a))$ and the induced representation on this space decomposes into a direct sum of irreducible representations with central character $\o$, all with finite multiplicities ([GP-S], [GGP-S], [Go1]). Such an irreducible subrepresentation is called {\bf cuspidal}. The sum $L_0^2(G(F)Z(\a)\bc G(\a))$ of the cuspidal representations is the {\bf cuspidal spectrum}.

The sum of all irreducible subrepresentations of $L^2(G(F)Z(\a)\bc G(\a))$ is  called the {\bf discrete spectrum} and it is denoted by $L^2_{disc}(G(F)Z(\a)\bc G(\a))$. It decomposes into a direct sum of irreducible representations with finite multiplicities. An irreducible subrepresentation of $L^2_{disc}(G(F)Z(\a)\bc G(\a))$ is called a {\bf discrete series}. The space $L^2_{disc}(G(F)Z(\a)\bc G(\a))$ admits an orthogonal decomposition
$$L^2_{disc}(G(F)Z(\a)\bc G(\a)) = L^2_0(G(F)Z(\a)\bc G(\a))\ \oplus\ L^2_{res} (G(F)Z(\a)\bc G(\a))$$
where $L^2_{res} (G(F)Z(\a)\bc G(\a))$ is stable by the action $R$ of $G(\a)$ and  decomposes into a sum of irreducible representations, all with finite multiplicities. It is called the {\bf residual spectrum}. Finally, one has orthogonal decompositions:
$$L^2(G(F)Z(\a)\bc G(\a))=$$
$$L_{disc}^2(G(F)Z(\a)\bc G(\a))\ \oplus\ L_{cont}^2(G(F)Z(\a)\bc G(\a))=$$
$$L_0^2(G(F)Z(\a)\bc G(\a))\ \oplus\  L_{res}^2(G(F)Z(\a)\bc G(\a))\ \oplus\  L_{cont}^2(G(F)Z(\a)\bc G(\a)),$$
where $L_{cont}^2(G(F)Z(\a)\bc G(\a))$ is a continuous sum of representations which has no irreducible quotient (and referred to as the {\bf continuous spectrum}).\\


If $(\pi,W)$ is a discrete series of $G(\a)$ ($\pi$ acts by right translations on the classes of functions in the subspace $W$ of $L^2(G(F)Z(\a)\bc G(\a))$), then for every $f\in \H(G(\a))$ we define $\pi(f):W\to W$ by $\pi(f)(\phi)=\phi'$, where, for all $x\in G(\a)$, 
$$\phi'(x)=\int_{G(\a)}f(g)\phi(xg)dg.$$
Because of the Lemma \ref{negligible} this is well defined, in the sense that, if $\psi$ is in the class of $\phi$, then $\psi'$ is in the class of $\phi'$.

To $\pi$ we may associate a set $\{\pi_v\}_{v\in V}$, such that 

- for all $v\in V$, $\pi_v$ is a unitary irreducible smooth representation of $G_v$, 

- for all but a finite number of $v\in V$, $\pi_v$ is unramified, and 

- for every simple tensor $f=\otimes_{v\in V}f_v\in \H(G(\a))$ one has 
$$\tr \pi (f)=\prod_{v\in V}\tr \pi_v(f_v).$$
I am not sure there is a place in the literature where a complete proof of this essential result is given. It may be recovered putting together results from  [H-C1] and [BJ] with  [Fl1], [GGP-S] and [KnVo] (see also [Bu], section 3.4).

Notice that, if $(\pi_v,E)$ is unramified and $f_v$ is the characteristic function of $K_v$, then $\pi_v(f_v)$ is a projection on the space of $V^{K_v}$. When $K_v$ is hyperspecial (that is the case for all but a finite number of $v\in V$), the dimension of $V^{K_v}$ is one and  $\tr\pi_v(f_v)=1$ (see Lemma \ref{spherical}).\\

Because of the independence of distribution characters for $G_v$, it is not hard to see that for every $v\in V$, the representation $\pi_v$ is  determined up to isomorphism by $\pi$.\\

The representation $\pi_v$ is called the {\bf local component of $\pi$ at $v$}. We say {\bf $\pi$ is ... at $v$} to say $\pi_v$ is ..., and {\bf $\pi$ is ... at almost every place} if, for almost all $v\in V$, $\pi_v$ is .... We know, for example, that $\pi$ is unramified at almost every place.

\section{Trace formula -- the compact quotient case}\label{compactquotient}

\def\fa{G(F)Z(\a)\bc G(\a)}
\def\lfa{L^2(G(F)Z(\a)\bc G({\a}))}
\def\lfazero{L_0^2(G(F) \bc G(\a))}

\def\oo{{\mathcal O}}
\def\cc{{\mathbb C}}

{\it In this section we assume the quotient $G(F)Z(\a)\bc G(\a)$ to be compact}. We make the general convention that, if $H$ is a group and $x\in H$, then $H_x$ is the centralizer of $x$ in $H$. The definition of the orbital integrals given below in this section is the same when the quotient is no longer compact, {\it if the integral defining it converges}.\\

Fix a smooth unitary character $\o$ of $Z(\a)$ trivial on $Z(F)$ and let $\lfa$ be like in the previous section the space of functions $\phi :G(\a) \to \cc$ such that $\phi(\g z g)=\o(z)\phi(g)$ for all 
$\g\in G(F)$, all $z\in Z(\a)$ and almost all $g\in G(\a)$, and 
$$\int_{\fa}|\phi(g)|^2dg<\infty.$$
Let $R$ be like in the previous section, i.e. the action of $G(\a)$ on $\lfa$ by right translations: $R(g)(\phi)(x)=\phi(xg)$.

For $f\in \H(G(\a))$ let $f_\o:G(\a)\to \cc$,
$$f_\o(g)=\int_{Z(\a)}\o(z)f(zg)\, dz.$$
That is well defined since the support of $f$ is compact.

One has 
$$R(f)(\phi)(x):= \int_{G(\a)} f(g)\phi(xg)\, dg =  \int_{Z(\a)\bc G(\a)} \int_{Z(\a)} f(zg)\phi(xzg)\, dz\, dg$$
$$=\  \int_{Z(\a)\bc G(\a)} \int_{Z(\a)} f(zg)\o(z)\phi(xg)\, dz\, dg\ =\ \int_{Z(\a)\bc G(\a)} f_\o(g)\phi(xg)\, dg.$$
\ \\

We want to interpret $R(f)$ as a kernel operator, and then apply the classical theorem on trace kernel operators which says:  {\it if $D$ is a measured space, if $k$ is a kernel operator in the space $W$ of complex functions on $D$ -- i.e. for all $\phi\in W$, for all $x\in D$,
$k(\phi)(x)=\int_D K(x;g)\phi(g)dg$, where $K:D\t D\to\cc$ is, by definition, the kernel -- then: if $k$ is trace class, if $K$ is continuous and   $K(g;g)$ is integrable over $D$, one has $\tr(k)=\int_D K(g;g)dg$.}\\
 
One has :
$$R(f)(\phi)(x)=\int_{Z(\a)\bc G(\a)} f_\o(g)\phi(xg)dg=\int_{Z(\a)\bc G(\a)} f_\o(x^{-1}g)\phi(g)dg$$
$$=\int_{\fa} \bigg({\sum_{\g\in Z(F)\bc G(F)}} f_\o(x^{-1}\g g)\phi(\g g)\bigg)\ dg$$
$$=\int_{\fa} \bigg({\sum_{\g\in Z(F)\bc G(F)}} f_\o(x^{-1}\g g)\bigg)\phi(g)\ dg.$$
If we set
$$K_f(x;g)={\sum_{\g\in Z(F)\bc G(F)}} f_\o(x^{-1}\g g),$$
then 
$$R(f)(\phi)(x)\ =\ \int_{\fa}\ K_f(x;g)\phi(g)\ dg.$$
Notice that $K_f$ is well defined: since $\o$ is trivial on $Z(F)$, $f_\o$ is $Z(F)$ invariant; also,
because $G(F)$ is a discrete subgroup of $G(\a)$ and  $f$ has compact support, given $x$ and $g$ the sum over $\g$ is finite.

Now $K_f$ is not strictly speaking the kernel of $R(f)$, because it is not well defined on $G(F)Z(\a)\bc G(\a)\t G(F)Z(\a)\bc G(\a)$
($K_f$ is defined on $G(F)\bc G(\a)\times G(F)\bc G(\a)$ and has the property that $K_f(z_1x;z_2g)=\o(z_1)\o^{-1}(z_2)K_f(x;g)$ for all $z_1,z_2\in Z(\a)$).

But one may chose a measurable fundamental domain $D$ for $G(F)Z(\a)\bc G(\a)$ (see for example [Bo] or [Go2]). Then, taking restriction to $D$ we may reinterpret the whole situation as:
$$R(f)(\phi)(x)\ =\ \int_D\ K_f(x;g)\phi(g)\ dg,$$
where $K_f$ is defined on $D\t D$. Now $K_f(g;g)$ is $Z(\a)$ invariant, and integrable over $D$ if and only if integrable over $G(F)Z(\a)\bc G(\a)$.

Recall we assumed $G(F)Z(\a)\bc G(\a)$ is compact (this is the case for example if $G(F)$ is a central division algebra of finite dimension over $F$; for a proof see [Booh]).
Then for all $f\in H(G(\a))$, $R(f)$ is a trace class operator, $K_f$ is continuous and $K_f(g;g)$ is integrable over $D$ and $G(F)Z(\a)\bc G(\a)$ (see for example [Booh] where all the details are carefully given). In the end we get:
$$\tr(R(f))\ =\ \int_D K_f(g;g)\, dg\ =\ \int_{\fa} K_f(g;g)\ dg$$
$$=\ \int_{\fa} \bigg({\sum_{   \g\in Z(F)\bc G(F)}} f_\o (g^{-1}   \g g)\bigg)\, dg.$$
\ \\

Let $X$ be the set of conjugacy classes in $Z(F)\bc G(F)$. Fix, for every $\oo\in X$, a representative $\g_o$ in $G(F)$ and a Haar measure on the centralizer $G_{\g_o}(\a)$.

For every conjugacy class $\oo\in X$, set
$$K_{f,\oo}(x;g)=\sum_{   \g\in\oo}f_\o(x^{-1}   \g g).$$ 
Then $K_f=\sum_{\oo}K_{f,\oo}$. For $\oo\in X$ one has
$$J_{\oo}(f)= \int_{\fa} K_{f,\oo}(g;g)dg=
\int_{\fa} \bigg(\sum_{\g\in\oo} f_\o (g^{-1} \g g)\bigg)dg =$$
$$ \int_{G_{\g_o}(F)Z(\a)\bc G(\a)} f_\o (g^{-1}\g_o g)\ dg \ =\ 
\vol(G_{\g_o}(F)Z(\a)\bc G_{\g_o}(\a))\int_{G_{\g_o}(\a)\bc G(\a)} f_\o (g^{-1}\g_o g)dg.$$ 

For $x\in G(\a)$ and $f\in \H(G(\a))$, we let
$$\Phi(f,x):=\int_{{G_x}(\a)\bc G(\a)} f (g^{-1}x g)dg$$
be {\bf the orbital integral of $f$ at $x$} and
$$\Phi(f_\o,x):=\int_{{G_x}(\a)\bc G(\a)} f_\o (g^{-1}x g)dg$$
be the {\bf orbital integral of $f_\o$ at $x$}. We have:
$$\Phi(f_\o,x)=\int_{Z(\a)}\Phi(f,zx)dz.$$
\ \\

We obtained the {\bf trace formula in the compact modulo the center case}:
$$\tr(R(f))\ =\ \sum_{\oo\in X} \vol(G_{\g_o}(F)Z(\a)\bc G_{\g_o}(\a))\ \Phi (f_\o,\g_o )\ =$$
$$\sum_{\oo\in X} \vol(G_{\g_o}(F)Z(\a)\bc G_{\g_o}(\a))\ \int_{Z(\a)} \Phi(f,z\g_o )dz.$$

The quotients are compact and so the volumes are finite. They depend on the choice of a measure on $G_{\g_o}$, but so does the orbital integrals and changing measures does not change the trace formula.

Hence, the trace of the operator $R(f)$ is a sum of orbital integrals of $f$. 
$\tr(R(f))$ is called {\bf the spectral side} of the trace formula, and the sum of orbital integrals is called {\bf the geometric side}. The second is theoretically more easy to compute. Still under the assumption $\fa$ is compact, the representation $R$ breaks into a discrete sum of unitary irreducible representations each of them appearing with finite multiplicity ([Booh] Corollary 2.2.2). All these representations are cuspidal (indeed, because we are in the compact quotient case, $G$ does not have proper parabolic subgroups). As these irreducible spaces are stable by  $R(f)$, $\tr(R(f))$ decomposes into a sum according to this decomposition. 

Let $f=\otimes_{v\in V}f_v$ be a simple tensor. We have seen that, if $\pi$ is an irreducible subrepresentation of $R$, then 
$$\tr\pi(f)=\prod_{v\in V}\tr \pi_v(f_v).$$ 
Also, one may prove that 
$$\Phi(f,x)=\prod_{v\in V} \Phi(f_v,x_v),$$ 
if $x=(x_v)_{v\in V}$, where
$$\Phi(f_v,x_v):=\int_{G_{v,x_v}\bc G_v} f_v(g_v^{-1}x_vg_v)dg_v,$$
with $G_{v,x_v}$ the centralizer of $x_v$ in $G_v$, is the {\bf orbital integral of $f_v$ in $x_v$}. (One has to interpret the quotient $G_x(\a)\bc G(\a)$ as the restricted product of the $G_{v,x_v}\bc G_v$, $v\in V$, with respect to the $K_v\cap G_{v,x_v}\bc G_v$, $v\in V_f$, and fix compatible measures on the quotients; this is left as an exercise). Notice that, when we wrote $x=(x_v)_{v\in V}$, we used the standard embedding of $G(F)$ in $G(\a)$ (section \ref{adelgroup}), which means that, for every $v$, one has $x_v=x$ (using the inclusion $G(F)\subset G_v$).\\

What happens when $\fa$ is no longer compact ? Consider the example of $GL_2({\q})$. 
If, for example, we want to compute $J_{\oo}$ for the conjugacy class 
$\oo$ of the element 
$\g=
\begin{pmatrix}
1&0\\
0&2
\end{pmatrix}
$ we find
$$ 
J_{\oo}(f)=\vol(G_{\g}(F)Z(\a) \bc G_{\g}(\a))\int_{G_{\g}(\a) \bc
G(\a)}f_\o (g^{-1}\g g)dg,$$
and, as $G_{\g}$ is the diagonal torus, one has

$$G_{\g}(F)Z(\a) \bc G_{\g}(\a)\simeq \ \q^\t\bc \a^\t$$
which has infinite volume (see the end of Section \ref{restprodmeas}). So $J_{\oo}(f)$ has no meaning.\\

\chapter{General linear groups}

\section{General linear groups over division algebras}\label{genlingroups}

Let $F$ be a number field and $V$ the set of places of $F$ with the usual notation $V_f$ and $V_\infty$ (Chapter 1, Section \ref{adelgroup}). 
Let $D$ be a central division algebra of finite dimension over $F$. Let $d^2$, $d\in\n^*$, be the dimension of $D$ over $F$ (it is always a square). Let $n\in\n^*$. 
Set $A:=M_{nd}(F)$ and $A':=M_n(D)$. For $v\in V$, let $A_v$, $A'_v$, be the completion of $A$, respectively of $A'$, at the place $v$. For any $v\in V$, one has $A_v\simeq M_{nd}(F_v)$ and $A'_v\simeq M_{nr_v}(D_v)$ where $r_v$ is a divisor of $d$ and $D_v$ is a central division algebra of dimension $d_v^2$ over $F_v$, with $d_v=d/r_v$\footnote{If $v\in V_\infty$, then $F_v$ is isomorphic to $\r$ or $\cc$ and $D_v$ is isomorphic with $\r$, $\cc$ or the quaternion algebra $\h$. When $v\in V_f$ every perfect square is the dimension of at least one division algebra over $F_v$.}. 
Notice that $D_v$ is not in general a completion of $D$. If $v\in V$, we say that $A'$ {\bf splits} at $v$ if $d_v=1$. The set of places $v$ of $F$ where $A'$ does not split is finite and we denote it by $S$.
So, if $v\in V\bc S$, one has $A_v\simeq M_{nd}(F_v)\simeq A'_v$. 
When $v\in V$, we fix once and for all isomorphisms $A_v\simeq M_{nd}(F_v)$ and $A'_v\simeq M_{nr_v}(D_v)$ which allow us to treat elements of $A_v$ and $A'_v$ as matrices, and to consider $A_v=A'_v$ for all $v\in V\bc S$. Thanks to the Skolem-Noether theorem, any algebra automorphisms of $A_v$ and $A'_v$ is inner ([Kn2], Corollary 2.42), and what follows does not depend in any essential way on the isomorphisms we fixed. A reference for this paragraph is [We].\\

Set $G:=A^\t=GL_{nd}(F)$ and $G':={A'}^\t=GL_n(D)$. For all $v\in V$, let $G_v= A_v^\t=GL_{nd}(F_v)$ and $G'_v={A'}_v^\t=GL_{nr_v}(D_v)$. 
We say $G'$ {\bf splits} at $v$ if $A'$ splits at $v$ (if $G'$ splits at $v$, then $G_v=G'_v$).\\ 

If $v\in V_f$, let $O_v$ and $O'_v$ be the ring of integers of $F_v$ and $D_v$ respectively. Set $K_v:=GL_n(O_v)$ and $K'_v:=GL_{nr_v}(O'_v)$ and fix Haar measures on $G_v$ and $G'_v$ such that the volume of $K_v$ and $K'_v$ is one. {\it For simplicity we will assume from now on that $V_\infty\cap S=\emptyset$, i.e. $A'$ splits at all infinite places.} For every $v\in V_\infty$, we fix a Haar measure on $G_v=G'_v$. 

One may see $G$ and $G'$ as algebraic groups defined over $F$ and denote by $G(\a)$ and $G'(\a)$ the ad\`ele associated groups (see Chapter 1, Section \ref{adelgroup}). Then $G(\a)$ (respectively $G'(\a)$) is the restricted product of the $G_v$, $v\in V$, with respect to the $K_v$, $v\in V_f$ (respectively of the $G'_v$, $v\in V$, with respect to the $K'_v$, $v\in V_f$).

The local Jacquet-Langlands correspondence is a correspondence between the representations of $G_v$ and $G'_v$. The global Jacquet-Langlands correspondence is a correspondence between the discrete series of the ad\`ele groups $G(\a)$ and $G'(\a)$.\\

\section{Regular semisimple elements, centralizers, orbital integrals, unitary induced representations}\label{orbint}

Let $v\in V$. {\it Only for this section, we denote by $A$ the algebra $M_n(D)$ for some $n\in\n^*$ and some central division algebra 
$D$ of dimension $d^2$ over $F_v$, and  set $G:=A^\t=GL_n(D)$. 
In particular, the definitions and results of this section apply to both $G_v$ and $G'_v$.}

There is a way of defining a {\bf characteristic polynomial} $P_g$ for any element $g$ of $A$, like in the case $D=F_v$. Actually, there are several equivalent ways, but I will use Bourbaki's definition -- [Bou3] section 17, Definition 1 -- where the characteristic polynomial is called reduced characteristic polynomial: for $g\in A$, if $P$ is the characteristic polynomial of the $F_v$-endomorphism of $A$ given by left multiplication with $g$, 
then $P_g$ is the unique unitary polynomial in $F_v[X]$ such that $P_g^{nd}=P$. Then the characteristic polynomial (is well defined and) has the following properties which will be freely used in the sequel:\\

\begin{prop}\label{serre}
{\rm (a)} There is an extension $E$ of $F_v$ such that $[E:F_v]=d$ and there is a $F_v$-algebras isomorphism $f: E\otimes D\simeq M_d(E)$. 

{\rm (b)} Let $\tilde{f}$ be the associated embedding of $A=M_n(D)$ into $M_{nd}(E)=M_n(M_d(E))$ using $A\to E\otimes A$, $x\mapsto 1\otimes x$. If $g\in A$, then the characteristic polynomial of $\tilde{f}(g)$ as a matrix in $M_{nd}(E)$ has coefficients in $F_v$ and equals $P_g$. 
\end{prop}
\ \\
{\bf Proof.} (a) By [Serre] XII.2 Proposition 2, there is an extension $E$ of $F_v$ such that $[E:F_v]=d$ and $D\otimes E\simeq M_d(E)$ as $F_v$-algebras. The equality $[E:F_v]=d$ is not part of the claim of the proposition in [Serre], but it is shown in the proof.

(b) [Bou3], 17.3, Corollary 1 to Proposition 4.\qed
\ \\

\begin{prop}\label{polcar} Let $g\in A$ and let $P_{min,g}$ be the minimal polynomial of $g$ over $F_v$. Then

{\rm (a)} $P_g$ is of degree $nd$,

{\rm (b)} $g\in A^\t$ if and only if $P_g(0)\neq 0$,

{\rm (c)}  $P_g(g)=0$, i.e. $P_{min,g}$ divides $P_g$.

{\rm (d)} $P_g$ divides $(P_{min,g})^{nd}$. In an algebraic closure of $F_v$, $P_g$ and $P_{min,g}$ have the same roots.

{\rm (e)} If $g'\in A$ is conjugate to $g$, then $P_{g'}=P_g$.
\end{prop}
\ \\
{\bf Proof.} (a) [Bou3], remark under the Definition 1 section 17.2.

(b) [Bou3] section 17.2, Proposition 3 combined with the Definition 2 of the reduced norm.

(c) Let $\tilde{f}$ be like in Proposition \ref{serre} (b). By Hamilton-Cayley, $P_g(\tilde{f}(g))=0$. So $P_g(g)=0$.

(d) With the notation from (c):

- If $Q$ is the minimal polynomial of $\tilde{f}(g)$ over $E$ then $P_g|Q^{nd}$ by standard matrix theory over fields. 

- Because $E$ is an extension of $F_v$,  $Q$ divides the minimal polynomial of $\tilde{f}(g)$ over $F_v$. 

- The minimal polynomial of $\tilde{f}(g)$ over $F_v$ is $P_{min,g}$ because $\tilde{f}$ is a morphism of $F_v$-algebras.\\ 
The rest follows using (c).

(e) Easy, using either the definition or Proposition \ref{serre}.\qed   
\ \\

If $(n_1,n_2,...,n_k)$ is a partition of $n$, we denote by $H(n_1,n_2,...,n_k)$ the subalgebra of $A$ made of diagonal matrices by blocks of size $n_1,n_2,...,n_k$. We set $L(n_1,n_2,...,n_k):=H(n_1,n_2,...,n_k)\cap G=H(n_1,n_2,...,n_k)^\t$ and call it the {\bf standard Levi subgroup} of $G$ associated to $(n_1,n_2,...,n_k)$. We denote by $(g_1,g_2,...,g_k)$ an element of $H(n_1,n_2,...,n_k)$, understood that $g_i\in M_{n_i}(D)$ for all $i\in\{1,2,...,k\}$. Then the characteristic polynomial of an element $(g_1,g_2,...,g_k)$ of $H(n_1,n_2,...,n_k)$ is the product of the characteristic polynomials of the elements $g_i\in M_{n_i}(D)$. 

An element $g\in G$ is called {\bf regular semisimple} if the characteristic polynomial $P_g$ of $g$ has simple roots (i.e. when decomposing $P_g$ in an algebraic closure of $F_v$, all roots appear with multiplicity one). Notice that $P_g$ is then the minimal polynomial of $g$ (by Proposition \ref{polcar} (c) and (d)). We denote by $G^{rs}$ the set of regular semisimple elements of $G$. A conjugacy class $c$ of $G$ is called {\bf regular semisimple} if it contains a regular semisimple element of $G$. Then all the elements of $c$ are regular semisimple and, if $g\in c$, then the elements of $c$ are exactly the elements of $G$ having a characteristic polynomial equal to $P_g$ (see the proof of Lemma 2.1 in [BaRo]).

If $P\in F_v[X]$ and $u\in \n^*$, we say $P$ is {\bf $u$-compatible} if the degree of any divisor of $P$ is divisible by $u$.

\begin{prop}\label{polcharbij}
{\rm (a)} If $g\in {G}^{rs}$ and $P_g$ is the characteristic polynomial of $g$, then $P_g$ is $d$-compatible.

{\rm (b)} Conversely, if $P\in F_v[X]$ is unitary, of degree $nd$, with simple roots, all non zero, and $d$-compatible, then $P$ is the characteristic polynomial of a regular semisimple element of $G$.
\end{prop}
\ \\
{\bf Proof.}   See [BaRo] Lemma 2.1.\qed 
\ \\

\begin{cor}\label{corpolcharbij}
Let $g\in G^{rs}$ and let $P_g$ be the characteristic polynomial of $g$. Let $P_g=\prod_{i=1}^kP_i$ be the decomposition of $P_g$ in a product of irreducible unitary polynomials in $F_v[X]$. Then, for every $i\in\{1,2,...,k\}$, $d$ divides the degree $\deg P_i$ of $P_i$. Moreover, if $\deg P_i=dn_i$, then $P_i$ is the characteristic polynomial of some $g_i\in GL_{n_i}(D)$, and $g$ is conjugate to the element $(g_1,g_2,...,g_k)$ of $L(n_1,n_2,...,n_k)$. 
\end{cor}
\ \\
{\bf Proof.} The fact that $d$ divides $\deg P_i$ follows from Proposition \ref{polcharbij} (a). The existence of the $g_i$ is a consequence of Proposition \ref{polcharbij} (b). Then $(g_1,g_2,...,g_k)$ and $g$ are conjugate because they have the same characteristic polynomial and $g$ is regular semisimple (see the proof of Lemma 2.1, [BaRo]).\qed
\ \\

\begin{lemme}\label{skolemnoether}
Let $g\in G^{rs}$. Let $U_g$ be the subalgebra of $A$ spanned by $g$ and $A_g$ be the centralizer of $g$ in $A$. Then

{\rm (a)} One has $U_g=A_g$ and $\dim_{F_v} U_g=nd$. Moreover, $U_g$ is a product of fields.

{\rm (b)} If $h\in G^{rs}$ is such that there is an isomorphism $i:U_h\simeq U_g$, then $i$ is conjugation with an element of $G$.
\end{lemme}
\ \\
{\bf Proof.}

(a) Let $g\in G^{rs}$ and $P_g=\prod_{i=1}^k P_i$ as in the corollary. Let $U_g$ be the $F_v$ subalgebra of $A$ spanned by $g$. Because $P_g$ is the minimal polynomial of $g$, $U_g\simeq F_v[X]/(P_g)$. The polynomials $P_i$, $i=1,2,...,k$, are pairwise relatively prime since $P_g$ has simple roots. By the Chinese lemma, $U_g\simeq \t_{i=1}^k K_i$ with $K_i=F_v[X]/(P_i)$ a field. Because a product of fields is a semisimple algebra, because  $\dim_{F_v}U_g=\deg(P_g)=nd$, by [Bou3] section 14.6 Proposition 3, $U_g$ is a maximal commutative sub-algebra of $A$. So $U_g$ must be the centralizer of $g$ (if $ag=ga$, then the subalgebra spanned by $U_g$ and $a$ is commutative).

(b) follows from (a) and the Corollary to Proposition 3, [Bou3] section 14.6.\qed
\ \\

Let $g\in G^{rs}$ as before. The centralizer $T_g$ of $g$ in $G$ is a maximal torus of $G$, and every maximal torus is the centralizer of one of its elements. One has 
$$T_g\ =\ U_g\cap G\ =\ U_g^\t\simeq(F_v[X]/(P_g))^\t\ \simeq\ K_1^\t\t K_2^\t\t ...\t K_n^\t.$$ 
If $g,g'\in G^{rs}$ are conjugate in $G$, then obviously the centralizers $T_g$ and $T_{g'}$ are conjugate in $G$.
We will fix measures on all the maximal tori of $G$, such that when $T$ and $T'$ are conjugate, the measure of $T'$ is the image of the measure of $T$ under the conjugation.  Let us be more precise:
let $T$ be a maximal torus of $G$ and let $dt$ be a Haar measure on $T$. Let $g\in G$. One may define a measure $dt_g$ on $gTg^{-1}$ in the following way:

- if $X\subset gTg^{-1}$, then $X$ is measurable if and only if $g^{-1}Xg$ is measurable in $T$ and then

- $\vol(X,dt_g):=\vol(g^{-1}Xg,dt)$.\\

\begin{lemme}\label{conjtori}
If $h\in G$ is such that $hTh^{-1}=gTg^{-1}$, then $dt_h=dt_g$.
\end{lemme}
\ \\
{\bf Proof.} Set $g^{-1}h=u$. Then $uTu^{-1}=T$. It is enough to show that $dt_u=dt$. 
Because $T$ is a maximal torus of $G$, it is known that the group $W$ of automorphisms of $T$ which are given by the conjugation with an element of $G$ is finite. For every $w\in W$, there is $\delta(w)\in \r_+^*$ such that $dt_w=\delta(w)dt$ (unicity of the Haar measure). It is easy to see that $\delta:W\to\r_+^*$ is a group morphism. As the only finite subgroup of $\r_+^*$ is $\{1\}$, it follows that $\delta$ is trivial and 
so $dt_w=dt$ for every $w\in W$. In particular, $dt_u=dt$.\qed
\ \\

We will simply say that $dt_g$ is {\bf the measure on $gTg^{-1}$ obtained from $dt$ by conjugation}, since it does not depend on the conjugation.
Now chose a set $S$ of representatives of conjugacy classes of maximal tori of $G$. For every $T\in S$, fix a Haar measure on $T$. If $T'$ is a maximal torus of $G$, there is a unique $T$ in $S$ such that there exists $g\in G$ such that $T'=gTg^{-1}$. On $T'$ we fix the measure obtained from $dt$ by conjugation. With this choice on every maximal torus of $G$, if two maximal tori are conjugate in $G$, their measures are obtained from one another by conjugation.\\

With this choice of measures, we define the orbital integrals. Recall that, if $f\in \H (G)$ and $g\in G^{rs}$, the {\bf orbital integral} of $f$ at  $g$ is by definition the complex number:
$$\Phi(f,g):=\int_{T_g\bc G} f(x^{-1}gx)\bar{dx},$$
where $\bar{dx}$ is the quotient measure. The function $g\mapsto \Phi(f,g)$ is stable by conjugation. 
\ \\

We gather now some facts concerning parabolic induction. We identify the standard Levi subgroup $L:=L(n_1,n_2,...,n_k)$ with $GL_{n_1}(D)\t GL_{n_2}(D)\t ...\t GL_{n_k}(D)$ and set $L^{rs}:=GL_{n_1}(D)^{rs}\t GL_{n_2}(D)^{rs}\t ...\t GL_{n_k}(D)^{rs}$.

Let $g=(g_1,g_2,...,g_k)\in G^{rs}\cap L$. Because the characteristic polynomial $P_g$ of $g$ has simple roots and $P_g=\prod_{i=1}^kP_{g_i}$, all the $P_{g_i}$ have simple roots and $g_i\in GL_{n_i}(D)^{rs}$ for all $i\in\{1,2,...,k\}$. So $G^{rs}\cap L\subset L^{rs}$. 
The converse is not true: if $g_i$ is regular semisimple in $GL_{n_i}(D)$ for all $i\in\{1,2,...,k\}$, $g$ is not regular semisimple in $G$ if the $P_{g_i}$ have common roots (hence common factors). However, we have the following:

\begin{prop}\label{essential}
{\rm (a)} $G^{rs}\cap L$ is a dense subset of $L^{rs}$, 

{\rm (b)} the centralizer of an element $g\in L^{rs}$ in $L$ equals the centralizer of a regular semisimple element of $G$ and is a maximal torus of $G$.
\end{prop}
\ \\
{\bf Proof.}  Let $g=(g_1,g_2,...,g_k)\in L^{rs}$ as before. Assume $g$ is regular semisimple in $L$ but not in $G$. We have the relation $P_g=\prod_{i=1}^kP_{g_i}$ of characteristic polynomials and some polynomials $P_{g_i}$ and $P_{g_j}$, $i\neq j$, have common factors. If $\varepsilon=(\l_1Id_{n_1}, \l_2 Id_{n_2},...,\l_k Id_{n_k})$ is an element of the center of $L$ with $\l_i\in F_v^\t$, then one may chose $\l_i$, as small as we want, such that the polynomials $P_{g_i+\l_i Id_{n_i}}$ are pairwise relatively prime, i.e. $g+\varepsilon\in L\cap G^{rs}$. This shows that $L\cap G^{rs}$ is dense in $L^{rs}$. Moreover, $\varepsilon$ being in the center of $L$, the centralizer of $g$ in $L$ equals the centralizer of $g+\varepsilon$ in $L$. Because $g+\varepsilon\in G^{rs}$, its centralizer in $G$ is included in $L$ (apply Lemma \ref{skolemnoether} (a)) and it is a maximal torus.\qed 
\ \\

In particular, as we fixed Haar measures on centralizers of regular semisimple elements of $G^{rs}$, we also fixed Haar measures on centralizers of regular semisimple elements of all standard Levi subgroups.\\
\ \\ 

We give some results about {\it unitary} representations {\it which are the only relevant representations for this paper}. All the representations are assumed to be admissible.


\begin{theo}\label{unitirred}
Let $P$ be a parabolic subgroup of $G$. Let $P=LN$ be a Levi decomposition of $P$. 
Let $\pi$ be a unitary irreducible representation of $L$. Then the induced representation $\ind_P^G\pi$ is unitarizable irreducible.
\end{theo} 
\ \\
{\bf Proof.} The unitarizability of $\ind_P^G\pi$ is known to follow easily from the definition of $\ind$ ([Re] IV.2.3, since $P\bc G$ is compact).
The irreducibility of $\ind_P^G\pi$ is a difficult result. It is proved in 
[Be1] when $D$ is a non archimedean field, [Se] when $D$ is non archimedean and not a field, [Bar] when $D$ is an archimedean field and [Vo] when $D$ is the quaternion algebra over $\r$ (see [BaRe] section 12).\qed
\ \\

\begin{prop}\label{pequalq}
Let $P$, $Q$ be  parabolic subgroups of $G$ with Levi decomposition $P=LN$ and $Q=MU$, where $N$ {\rm (}respectively $U${\rm )} is the unipotent radical of $P$ {\rm (}respectively of $Q${\rm )}. Let $\pi$ be a representation of $L$. Assume there exists $g\in G$ such that $M=gLg^{-1}$, and let $\pi^g$ be the representation of $M$ given by $\pi^g(x)=\pi(g^{-1}xg)$. If $\ind_P^G\pi$ is irreducible then $\ind_Q^G\pi^g$ and $\ind_P^G\pi$ are isomorphic. 
\end{prop} 
\ \\
{\bf Proof.} First assume $Q=gPg^{-1}$. Then, for $x\in G$, $(\ind_Q^G\pi^g)(x)\ =\ (\ind_P^G\pi)(g^{-1}xg)$, so the representations are isomorphic.
That brings us to the case when $M=L$.
When $F_v$ is non archimedean, one may use the formula of the induced character from [Cl] (Proposition 3) to prove that the character function of the induced representation does not depend of the choice of the parabolic subgroup. The same formula works when $F_v$ is archimedean ([Kn1], Proposition 10.18; appears as Theorem 10 in Delorme's survey [Del])).\qed
\ \\

In particular, if $P=LN$ and $\pi$ is a unitary irreducible representation of $L$, then $\ind_P^G\pi$ does not depend on $P$ and we simply denote it by $\ind_L^G\pi$.\\

\section{Transfer of conjugacy classes}\label{transconjclass}

Let $v\in V_f$. We go back to the notation $G_v$, $G'_v$. 
We write $g\lra g'$ if $g\in G_v^{rs}$, $g'\in {G'}_v^{rs}$ and we have equality $P_g=P_{g'}$ of characteristic polynomials. If $g\in G^{rs}$ and $g'\in {G'}^{rs}$, we say $g$ {\bf corresponds to} $g'$ if $g\lra g'$. Then any element in the conjugacy class of $g$ corresponds to any element in the conjugacy class of $g'$. We obtain (Proposition \ref{polcharbij}) an injective map from the set of regular semisimple conjugacy classes of $G'$ into the set of regular semisimple conjugacy classes of $G$. Moreover, if $g\in G^{rs}_v$, we have an equivalence between

(i) there exists $g'\in {G'}_v^{rs}$ such that $g\lra g'$ and 

(ii) the characteristic polynomial of $g$ splits into a product of irreducible polynomials (over $F_v$) all of which degrees are divisible by $d_v$.\\
We denote $G_v^{rsd_v}$ the set of elements of $G^{rs}_v$ satisfying this condition.\\

\section{Transfer of centralizers and measures}\label{transtori}

Let $v\in V_f$. {\it On maximal tori of $G_v$ we fix measures like before, such that when two maximal tori are conjugate the Haar measures are conjugate}. On maximal tori of $G'_v$ we will fix Haar measures associated to the ones fixed on maximal tori of $G_v$. We will show that, when two maximal tori of $G'_v$ are conjugate, the Haar measures we fixed are conjugate. 

Let $g\lra g'$. Then $P_g=P_{g'}$  and $U_g\simeq F_v[X]/(P_g)\simeq U_{g'}$. We fix an algebra isomorphism $i:U_g\to U_{g'}$. This gives rise to a group isomorphism between $T_g=U_g^\t$ and $T_{g'}=U_{g'}^\t$ which we denote again by $i$. Using $i$, we transfer the measure $dt$ of $T_g$ to a measure $dt'$ on $T_{g'}$ (if $X\subset T_{g'}$, then $X$ is measurable if and only if $i^{-1}(X)$ is measurable and then $\vol(X):=\vol(i^{-1}(X))$). Because of the Lemma \ref{skolemnoether} (b), if we chose a different isomorphism $\theta:U_g\simeq U_{g'}$, there exists $t\in G'_v$ such that $\theta\circ i^{-1}$ is conjugation by $t$. Then, by Lemma \ref{conjtori}, the measure we obtain on $T_{g'}$ using $\theta$ is the same we obtained using $i$. A priori this measure depends on the choice of the couple $g,g'$, but the following argument will show that this is not the case. Assume $h'\in {G'}^{rs}$ is such that $T_{h'}$ is conjugate to $T_{g'}$ and let $h\in G^{rs}$ such that $h\lra h'$. Put on $T_{h'}$ the measure obtained from the one on $T_h$ like before. Now $U_{h'}$ is conjugate (hence isomorphic) to $U_{g'}$ (exercise) and so $U_h$ is isomorphic to $U_g$. By Lemma \ref{skolemnoether} (b), $U_h$ and $U_g$ are conjugate, so $T_h$ and $T_g$ are conjugate. The measures on $T_h$ and $T_g$ are then conjugate. The independence on $i$ of our construction shows (composing $i$ with conjugations) that the measures on $T_{h'}$ and $T_{g'}$ are conjugate.

We say we obtained the measure on $T_{g'}$ {\bf transferring from $G_v$}.

\section{Transfer of functions}\label{transfct}

Let $v\in V_f$. In the previous section we fixed measures on the maximal tori of $G_v$ such that if two maximal tori are conjugate the measures are conjugate, then we transferred the measures on the maximal tori of $G'_v$ where we showed that the same holds.

If $f\in \H(G_v)$ and $f'\in \H(G'_v)$ we say $f$ {\bf corresponds to} $f'$ and write $f\lra f'$ if 

- the support of $f$ is included in $G_v^{rs}$ and the support of $f'$ is included in ${G'_v}^{rs}$,

- $\Phi(f,g)=\Phi(f',g')$ when $g\lra g'$, and 

- $\Phi(f,g)=0$ if $g\in G_v^{rs}$ does not correspond to any element of $G'_v$ (i.e. $g\in G_v^{rs}\ \bc \ G_v^{rsd_v}$).\\
This is well defined since orbital integrals are stable by conjugation (as the relation $g\lra g'$). The problem is : are there corresponding functions? The answer is yes:

\begin{prop}
If $f'\in\H(G'_v)$ has support included in ${G'_v}^{rs}$, there exists $f\in\H(G_v)$ such that $f\lra f'$. Moreover, one may chose $f$ to have support in $G_v^{rsd_v}$.

If $f\in \H(G_v)$ has support included in $G_v^{rs}$ and is such that the orbital integral of $f$ is zero on $G_v^{rs}\ \bc \ G^{rsd_v}_{v}$, then there exists $f'\in\H(G'_v)$ such that $f\lra f'$.
\end{prop}
\ \\

The proposition is known to be a consequence of the submersion principle of Harish-Chandra [H-C2] (see [BaRo], Lemme 2.3 and Prop.2.4 for the proof).\\

\section{The local Jacquet-Langlands correspondence}\label{localjl}

Let $v\in V_f$. All representations are supposed admissible.\\

\begin{theo}\label{jllocal} {\bf (The local Jacquet-Langlands correspondence for unitary representations.)}
Let $\pi$ be a unitary irreducible representation of $G_v$. Then $\pi$ is in one of the following situations:

- $\tr\pi(f)=0$ whenever $f\lra f'$ for some $f'\in \H(G'_v)$,

- there is a unique up to isomorphism representation $\pi'$ of $G'_v$ such that there is a sign $\varepsilon_\pi\in \{-1,1\}$
such that
$$\tr\pi(f)=\varepsilon_\pi \tr\pi'(f')$$
whenever $f\lra f'$ for some $f'\in \H(G'_v)$. Moreover, $\pi'$ is unitary and irreducible.
\end{theo}
\ \\
{\bf Proof.} It is Theorem 5.2 in [BHLS].\qed
\ \\

When $\pi$ is in the first situation we will say $\pi$ {\bf corresponds to zero}, when $\pi$ is in the second situation, $\pi$ {\bf corresponds to} $\pi'$.\\

In the theorem, one may replace the hypothesis {\it $\pi$ unitary} by {\it $\pi$ square integrable} (resp. {\it tempered}, resp. {\it elliptic}, resp. {\it ladder}) and then the conclusion by {\it $\pi'$ square integrable} [DKV] (resp. {\it tempered} [DKV], resp. {\it elliptic} [Ba2], resp. {\it ladder} [BLM]). When $\pi$ is square integrable, the first situation never holds and the sign depends only on $G'_v$. This version of the correspondence (i.e. for square integrable representations) is the original one and was proved in [DKV] after long efforts ([JL], [Ro1]). All the other versions crucially use this original result in their proofs. In [Ba2] it is proved that any smooth irreducible representation $\pi$ will be in one of the two situations, if we let $\pi'$ be not a representation but a virtual sum of representations. Deng showed [De] that there are cases when the sum contains at least two different terms (i.e. irreducible does not correspond to irreducible). The version I gave here, Theorem \ref{jllocal}, is the good one for the statement and proof of the global correspondence. Indeed, the local components of discrete series are always unitary (but not necessarily tempered).\\
\ \\

\section{Corresponding Levi, parabolic induction and correspondence}\label{levicorresp}

Let $v\in V_f$. If $L(n_1,n_2,...,n_k)$ is a standard Levi subgroup of $G'_v$ we say that it {\bf corresponds} to the standard Levi subgroup $L(d_vn_1,d_vn_2,...,d_vn_k)$ of $G_v$. If $L(m_1,m_2,...,m_k)$ is a standard Levi subgroup of $G_v$, then we say it {\bf corresponds to nothing} if there is at least one $m_i$ which is not divisible by $d_v$, and we say it {\bf corresponds to} $L(m_1/d_v,m_2/d_v,...,m_k/d_v)$ if $d_v$ divides $m_i$ for all $i\in \{1,2,...,k\}$.

Because the standard Levi subgroups are product of linear groups, it is easy to extend the Jacquet-Langlands correspondence to standard Levi subgroups $L$ and $L'$ such that $L$ corresponds to $L'$. As induced representations from unitary irreducible are unitary irreducible (Theorem \ref{unitirred}), it is important to know how parabolic induction interact with the local Jacquet-Langlands correspondence. Beyond being a natural question it is needed for the proof of the global correspondence.

\begin{prop}\label{inductionjl}
Let $L$ be a standard Levi subgroup of $G_v$. Let $\pi$ be a unitary irreducible representation of $L$. Then

{\rm (a)} If $L$ corresponds to nothing or $\pi$ corresponds to zero, then $\ind_L^{G_v}\pi$ corresponds to zero.

{\rm (b)} If $L$ corresponds to $L'$ and $\pi$ corresponds to $\pi'$, then $\ind_L^{G_v}\pi$ corresponds to $\ind_{L'}^{G'_v}\pi'$.
\end{prop}
\ \\
{\bf Proof.} Let $P$ be a parabolic subgroup having $L$ as a Levi subgroup. If $f\in \H(G_v)$, then one may define the constant term $f^P$ of $f$ along $P$ which is an element of $\H(L)$ such that 

(i) if $g\in L\cap G_v^{rs}$, then $D_{G/L}^{-1/2}(g)\ \Phi(f^P,g)=\Phi(f,g)$

(ii) $\tr\ \ind_L^{G_v}\pi(f)=\tr\ \pi(f^P)$.\\
See [Lau] Prop.4.3.11 and Lemma 7.5.7 for the definition of the function $D_{G/L}$, of $f^P$ (4.3.8) and proofs. The proofs of Laumon are in positive characteristic but work the same here. See also [FK], chapter 7, page 62.

Assume $L$ corresponds to $L'$. Let $P$ be the parabolic subgroup of $G_v$ containing $L$ and the upper triangular matrices; let $P'$ be the parabolic subgroup of $G'_v$ containing $L'$ and the upper triangular matrices. Let $N$ (resp. $N'$) be the unipotent radical of $P$ (resp. of $P'$). Then, if $l\in L$ and $n\in N$, the characteristic polynomial of $ln$ is the characteristic polynomial of $l$. That follows from the classical definition of the characteristic polynomial as a determinant (over a commutative field) and the fact that the determinant of a matrix in $P$ depends only on the coefficients of the matrix in the diagonal blocs corresponding to $L$. Then the result is true also for $l\in L'$ and $n\in N'$: it follows from Proposition \ref{serre} and the previous classical case. So, if $f\in \H(G_v)$ (resp. $f'\in \H(G'_v)$) has support included in $G_v^{rs}$ (resp. in ${G'}_v^{rs}$), then $f^P$ (resp. $f'^{P'}$) has support included in $L\cap G_v^{rs}$ (resp. $L'\cap {G'}_v^{rs}$). 

\begin{lemme}\label{lelongparab}
If $L$ corresponds to $L'$, if $f\lra f'$, then $f^P\lra {f'}^{P'}$.
\end{lemme}
\ \\
{\bf Proof.}  Because $f^P$ has support included in $G_v^{rs}\cap L$ and this set is stable by conjugation with elements of $L$, the orbital integral $\Phi(f^P,g)$ is zero if $g\in L^{rs}\bc G_v^{rs}$. For the same reasons, $\Phi({f'}^{P'},g')=0$ if $g'\in {L'}^{rs}\bc {G'_v}^{rs}$.

Now, if $g\in L^{rs}$ and $g'\in {L'}^{rs}$ are such that $g\lra g'$, then one has $g\in G^{rs}$ if and only if $g'\in {G'}^{rs}$.
If $g\in L\cap G_v^{rs}$ and $g'\in {L'}\cap {G'_v}^{rs}$ and $g\lra g'$, we have $\Phi(f^P,g)=\Phi({f'}^{P'},g')$ by (i) because $f\lra f'$ and $D_{G/L}(g)=D_{G'/L'}(g')$.

If $g\in L\cap G_v^{rs}$ does not correspond to any $g'\in L'$, then $g$ does not correspond to any $g'\in {G'}_v^{rs}$, so
 $\Phi(f^P,g)=0$ by (i).\qed 
\ \\

(a) First, assume that $L$ corresponds to nothing, i.e. among the sizes of the blocks of $L$ there is at least one not divisible by $d_v$. If $f\lra f'$, then the support of $\Phi(f,\ )$ is included in $G_v^{rsd_v}$. But $G_v^{rsd_v}$ has no intersection with $L$ because the characteristic polynomial of an element of $L$ is a product of polynomials having degrees equal to the sizes of blocks of $L$. 
So the orbital integral of $f^P$ is zero on $L\cap G_v^{rs}$ by (i). Because the support of $f^P$ is included in  $L\cap G_v^{rs}$ which is stable by conjugation with elements of $L$, the orbital integrals of $f^P$ are zero on the whole $L^{rs}$. Then $\tr\pi(f^P)$ is zero (the trace of a representation is zero an a function with everywhere null orbital integrals, [DKV] Theorem 2.f). So $\tr\ \ind_L^{G_v}\pi(f)=0$ by (ii).\\

Assume now $L$ corresponds to $L'$ and $\pi$ corresponds to zero. If $f\lra f'$, then, by the Lemma \ref{lelongparab}, $f^P\lra{f'}^{P'}$. But $\tr\ \pi(f^P)=0$ because $\pi$ corresponds to zero. By (ii),
$$\tr\ \ind_L^{G_v}\pi(f)\ =\ \tr \pi(f^P)\ =\ 0.$$
\ \\

(b) We saw that $f\lra f'$ implies $f^P\lra {f'}^{P'}$. Then, by (ii), 
$$\tr\ \ind_L^{G_v}\pi(f)\ =\ \tr\ \pi(f^P)\ =\ \tr\ \pi'({f'}^{P'})\ =\ \tr\ \ind_{L'}^{G'_v}\pi'(f').$$
\ \\

\section{The global Jacquet-Langlands correspondence}\label{globaljl}

Let $F$ be our number field and $D$ a central division algebra of dimension $d^2$ over $F$ as in Section \ref{genlingroups}. Let $S$ be the (finite) set of places of $F$ where $D$ does not split. Recall we assume here that $D$ splits at infinite places, i.e. $S\cap V_\infty=\emptyset$. We do that in order to avoid developing the harmonic analysis on real groups, for example the transfer of orbital integrals. This restriction is also made in [Ba1]. The global Jacquet-Langlands is proved without assuming $D$ splits at infinite places in [BaRe].

If $k\in\n^*$, let $G_k:=GL_k(F)$ and $G'_k:=GL_k(D)$. Then $G'_{k,v} = G_{kd,v}$ for all $v\in V\bc S$. Let $DS_k$ be the set of discrete series of $G_k(\a)$ and $DS'_k$ the set of discrete series of $G'_k(\a)$. 
The global Jacquet-Langlands correspondence is the following theorem:\\

\begin{theo}\label{jlglobal} {\bf (The global Jacquet-Langlands correspondence.)}
There exists an injective map ${\bf G}:DS'_k\to DS_{kd}$, such that, for every $\pi'\in DS'_k$, if $\pi:={\bf G}(\pi')$, then 

- $\pi_v\simeq \pi'_v$ if $v\in V\bc S$ and 

- $\pi_v$ corresponds to $\pi'_v$ by the local Jacquet-Langlands correspondence if $v\in S$.
\end{theo}
\ \\

\section{The six steps of the proof}\label{sixstepproof}

I will give a fake sketch, step by step, of the proof of Theorems \ref{jllocal} and \ref{jlglobal}. The steps are more or less the same as in the proof of any correspondence of representations, using the trace formula. {\bf The steps as described here in a first time are wishes, in the sense that things are not as simple as pretended}. For example, the groups $G_k(\a)$ and $G'_k(\a)$ do not have in general a trace formula as simple as in Step 2. However, these simplistic wishes have been the starting point of so many developments in the field. {\bf I will discuss every step after, explaining what is the real situation.} To simplify notation, in this section we fix $n$ and set $G=G_{nd}$ and $G'=G'_n$.\\
\ \\
\ \\
{\bf Step 1.} (State the correspondence. This requires definitions and proofs of transfer of conjugacy classes, of centralizers and their measures, of functions (i.e. of  orbital integrals).) 

Define a natural injective map $\mathcal{I}$ from the set of conjugacy classes of $G'_v$ into the set of conjugacy classes of $G_v$. 
Write $g\lra g'$ if the image of the class of $g'$ through $\mathcal{I}$ is the class of $g$. 

If $g\lra g'$, define in a coherent way measures on the centralizers of $g$ and $g'$. 

Write $f\lra f'$ if $f\in\H(G_v)$, $f'\in\H(G'_v)$ and 

- $\Phi(f,g)=\Phi(f',g')$ whenever $g\lra g'$ and

- $\Phi(f,\ )$ is zero on classes which are not in the image of $\mathcal{I}$.\\ 
Prove that for every $f'\in \H(G'_v)$ there is $f\in \H(G_v)$ such that $f\lra f'$.

Using that, state the correspondence, like for example in Theorems \ref{jllocal} and \ref{jlglobal}.
\ \\
Those things have been explained for general linear groups in the previous part of this chapter.\\
\ \\
\ \\
{\bf Step 2.} (The trace formula.) 

If $Z$ is the center of $G$ and $Z'$ is the center of $G'$, we have canonical isomorphisms $Z(\a)\simeq \a^\t\simeq Z'(\a)$ and we identify $Z'(\a)$ to $Z(\a)$. Then $Z(\a)$ is the restricted product of the $Z_v$, $v\in V$, with respect to the $Z_v\cap K_v$, $v\in V_f$.
Fix a smooth unitary character $\o$ of the center $Z(\a)$ of $G(\a)$ and $G'(\a)$. Let $R$ (resp. $R'$) be the representation of $G(\a)$ (resp. $G'(\a)$) by right translations in the space $L^2(G(F)Z(\a)\bc G(\a))$ (resp. $L^2(G'(F)Z(\a)\bc G'(\a))$) like in Chapter 1, Section \ref{discreteseries}.

Prove trace formulae:
$$\tr(R(f))=\sum_{\oo\in X} \vol(G_{\g_o}(F)Z(\a)\bc G_{\g_o}(\a))\int_{Z(\a)}\Phi(f,z\g_o)\ dz$$
and
$$\tr(R'(f'))=\sum_{\oo'\in X'} \vol(G'_{\g_{o'}}(F)Z'(\a)\bc G'_{\g_{o'}}(\a))\int_{Z'(\a)}\Phi(f',z\g_{o'})\ dz$$
for $G$ and $G'$, as in the compact quotient case, Chapter 1, Section \ref{compactquotient}.

Recall $X$ (resp. $X'$) is the set of conjugacy classes in $Z(F)\bc G(F)$ (resp. $Z(F)\bc G'(F)$) and for each element $\oo$ of $X$ (resp. $\oo'$ of $X'$) we chose an element $\g_o\in G(F)$ such that its class (modulo $Z(F)$) is in $\oo$ (resp. $\g_{o'}\in G'(F)$ such that its class is in $\oo'$).\\
\ \\
\ \\
{\bf Step 3.} (Comparison of the trace formulae.) 

Let $S$ be the set of places where $G'$ does not split. Then $G_v=G'_v=GL_{nd}(F_v)$ for $v\notin S$. If $f\in\H(G(\a))$ and $f'\in\H(G'(\a))$, write $f\lra f'$ if both are tensor products of local functions $f=\otimes_{v\in V} f_v$ and $f'=\otimes_{v\in V} f'_v$ such that 

- $f_v\lra f'_v$ for all $v\in S$, as defined at Step 1,

- $f_v=f'_v$ for all $v\notin S$.\\ 
Step 3 is to show the equality of the geometric sides of the trace formulae for a couple $f,f'$ such that $f\lra f'$:
$$\sum_{\oo\in X} \vol(G_{\g_o}(F)Z(\a)\bc G_{\g_o}(\a))\int_{Z(\a)}\Phi(f,z\g_o)\ dz$$
$$=\ \sum_{\oo'\in X'} \vol(G'_{\g_{o'}}(F)Z'(\a)\bc G'_{\g_{o'}}(\a))\int_{Z'(\a)}\Phi(f,z\g_{o'})\ dz$$
and so to get the equality $\tr(R(f))=\tr(R'(f'))$ of the spectral sides. 

The equality of the geometric sides will be proved term by term.
There is a natural injective map $J$ from the set of conjugacy classes of $G'(F)$ into the set of conjugacy classes of $G(F)$ preserving the characteristic polynomial. It induces an injective map from $X'$ to $X$, denoted by $\mathcal J$.\\

{\it We will show that the term indexed by $\oo'$ on the right hand side equals the term indexed by ${\mathcal{J}}(\oo')$ on the left hand side, and the other terms on the left hand side, indexed by orbits which are not in the image of ${\mathcal{J}}$, are zero.}\\

We will assume that the representative $\g_o$ of $\mathcal{J}(\oo')$ is in $J(\g_{o'})$ (that was already true up to multiplication with an element of the center), i.e. its characteristic polynomial $P_{\g_o}$ equals the characteristic polynomial $P_{\g_{o'}}$ of $\g_{o'}$.
Then the local components of $\g_o$ and $\g_{o'}$ verify $\g_{o,v}\lra \g'_{o,v}$ for all $v\in V$ (because $P_{\g_{o,v}}= P_{\g_o}=P_{\g_{o'}}=P_{\g_{o',v}}$) and $z_v\g_{o,v}\lra z_v\g'_{o,v}$ for all $z_v\in Z_v$.\\

Assume then measures on $G_{\g_o}(\a)$ and $G'_{\g_{o'}}(\a)$ are defined by the (restricted) product of local measures on local centralizers, of $\g_{o,v}$ in $G_v$ and $\g_{o',v}$ in $G'_v$. The measures on these local centralizers correspond to each other by the construction of Step 1. 
Then, when $\oo={\mathcal{J}}(\oo')$, one has 
$$\vol(G'_{\g_{o'}}(F)Z'(\a)\bc G'_{\g_{o'}}(\a))=
\vol(G_{\g_o}(F)Z(\a)\bc G_{\g_o}(\a)).$$
Also,
$$\Phi(f,z\g_o)=\prod_{v\in V} \Phi(f_v,z_v\g_{o,v})=\prod_{v\in V} \Phi(f'_v,z_v\g_{o',v})=\Phi(f,z\g_{o'}).$$

Now, when $\oo$ is a conjugacy class which is not in the image of $\mathcal{J}$, the corresponding term in the geometric side of the trace formula for $G(\a)$ is zero. Indeed, there is a place $w$ where $\g_{o,w}$ does not correspond to any element of $G'_w$. Then $z_w\g_{o,w}$, $z_w\in Z_w$, does not correspond to any element of $G'_w$.  Hence
$$\Phi(f,z\g_o)=\prod_{v\in V} \Phi(f_v,z_v\g_{o,v})=0,$$
because $\Phi(f_w,z_w\g_{o,w})=0$ by definition of $f_w\lra f'_w$. 

We proved the equality of the geometric sides of the trace formulae for $G$ and $G'$. So
$$\tr(R(f))=\tr(R'(f')).$$
\ \\
\ \\
{\bf Step 4.} (Separate discrete series with equal \ ${ }^S$-part.)

We write $R$ as a sum of irreducible non isomorphic representations $\pi_i$, $i\in I$, with finite multiplicity $m_i$ and the representation $R'$ as a sum of irreducible non isomorphic representations $\pi'_j$, $j\in J$, with finite multiplicities $m'_j$. By Step 3 we know that $\tr(R(f))=\tr(R'(f'))$, so 
$$\sum_{i\in I} m_i\tr\pi_i(f)\ =\ \sum_{j\in J} m'_j\tr\pi'_j(f')$$
if $f\lra f'$.

Recall $S$ is the set of places where $G'$ does not split, so $G(\a)=G_S\t G^S$ and $G'(\a)=G'_S\t G^S$, as $G'_v=G_v=GL_{nd}(F_v)$ for $v\notin S$. Every irreducible representation $\pi$ of $G(\a)$ may be written as 
$\pi=\pi_S\otimes \pi^S$ where $\pi_S$ is an irreducible representation of $G_S$ (the ${ }_S$-part of $\pi$) and $\pi^S$ is an irreducible representation of $G^S$ (the ${ }^S$-part of $\pi$). Every irreducible representation $\pi'$ of $G'(\a)$ may be written 
$\pi'=\pi'_S\otimes {\pi'}^S$ where $\pi'_S$ is an irreducible representation of $G'_S$ (the ${ }_S$-part of $\pi'$) and ${\pi'}^S$ is an irreducible representation of $G^S$ (the ${ }^S$-part of $\pi'$).

By linear independence of characters on the group $G^S$ (for $v\notin S$, we have $f_v=f'_v\in \H(G_v)=\H(G'_v)$), we may separate the equation  in equalities indexed by $t$, where $t$ runs over the irreducible representations of $G^S$, keeping on left and right only representations which have ${}^S$-part isomorphic to $t$:
$$\sum_{i\in I_t} m_i\tr\pi_i(f)\ =\ \sum_{j\in J_t} m'_j\tr\pi'_j(f')$$
where $I_t$ (resp. $J_t$) is the set of representations $\pi_i$, $i\in I$ (resp. $\pi'_j$, $j\in J$), such that $\pi_i^S\simeq t$ (resp. ${\pi'_j}^S\simeq t$). If we write $f=f_S\otimes f^S$ and $f'=f'_S\otimes f^S$ with the obvious notation, we have:
$$\sum_{i\in I_t} m_i\, \tr\pi_{i,S}(f_S)\, \tr(t(f^S))\ =\ \sum_{j\in J_t} m'_j\, \tr\pi'_{j,S}(f'_S)\, \tr(t(f^S)).$$
Then we may consider a function $f^S$ such that $\tr(t(f^S))\neq 0$, then we obtain
$$\sum_{i\in I_t} m_i\, \tr\pi_{i,S}(f_S)\ =\ \sum_{j\in J_t} m'_j\, \tr\pi'_{j,S}(f'_S).$$
\ \\
\ \\
{\bf Step 5.} (Prove a global correspondence.)

We prove strong multiplicity one theorems: {\it all the multiplicities $m_i$ and $m'_j$ are $1$, and given an irreducible representation $t$ of $G^S$, the cardinality of the sets $I_t$ and $J_t$ are at most one; i.e. if two discrete series are isomorphic at almost every place, then they are equal}. 

We fix a discrete series $\Pi'$ of $G'(\a)$ and set $t:={\Pi'}^S$.
Then the right side of the equality indexed by $t$ contains only one element, namely $\Pi'$. There is at least one function $f'\in\H(G'(\a))$ such that $\tr\Pi'(f')\neq 0$. This shows that the left side is not identically zero, so $I_t\neq\emptyset$. Then $I_t$ contains one element $\Pi$ (which is a discrete series of $G(\a)$) and we have:
$$\tr\Pi_S(f_S)=\tr\Pi'_S(f'_S)$$
if $f\lra f'$. We obtained a {\bf global Jacquet-Langlands correspondence}, that is an injective map {\bf G} from the set of discrete series of $G'(\a)$ into the set of discrete series of $G(\a)$ such that, if ${\bf G}(\Pi')=\Pi$, then $\tr\Pi(f)=\tr\Pi'(f')$ if $f\lra f'$.\\
\ \\
\ \\
{\bf Step 6.} (Prove the local correspondence.)

Only now, that we proved a global Jacquet-Langlands correspondence, we may come back to the local one and prove it. Fix a place $v_0\in V$. First, fix a function $f'$ such that $\tr\Pi'(f')\neq 0$, and then $f$ such that $f\lra f'$. The equality $\tr\Pi_S(f_S)=\tr\Pi'_S(f'_S)$ if $f_S\lra f'_S$ of Step 5 reads:
$$\prod_{v\in S}\tr\Pi_v(f_v)=\prod_{v\in S}\tr\Pi'_v(f'_v)$$
if $f_v\lra f'_v$ for all $v\in S$.
Then varying only $f'_{v_0}$ and $f_{v_0}$, such that $f_{v_0}\lra f'_{v_0}$, we obtain that there is a constant $\l\in \cc^*$ such that $\tr\pi_{v_0}(f_{v_0})=\lambda\tr\pi'_{v_0}(f'_{v_0})$ if $f_{v_0}\lra f'_{v_0}$. Finally, one may prove that $\l=\pm 1$. This is the local correspondence.\\
\ \\
\ \\
{\bf The truth.} Let us now explain what are the difficulties in accomplishing the steps 1 to 6. Keep in mind $G(\a)$ is not in the compact quotient case unless $n=d=1$, and $G'(\a)$ is not in the compact quotient case unless $n=1$.\\ 
\ \\
{\bf Step 1?} That has been already explained in Section \ref{transconjclass}, \ref{transtori} and \ref{transfct}. Notice however that we defined transfer only for regular semisimple elements. For proving correspondences, this transfer seems to be enough, and transfer of non regular semisimple classes may be avoided. The transfer is more complicated for other groups or other correspondences (see the first sections of [Ro2] for unitary groups and of [AC] for base change).\\
\ \\
\ \\
{\bf Step 2?} Such a formula may not be proved unless $G(\a)$ and $G'(\a)$ are in the compact quotient case (so, in general, not for our groups). The general trace formulae are much more complicated. 

However, the simple trace formula of Kazhdan and Deligne (see [DKV] or [He] 4.9) has this form, but applies only for a particular class of functions $f$ and $f'$ (they have to be "cuspidal" at one place $v_1$ and "supported in the regular elliptic set" at another place $v_2$). This simple trace formula applies to any reductive group, but has the inconvenient that the class of functions is not very large, in particular the trace of global representations on these functions is zero for discrete series which are not cuspidal (and even cuspidal at one place), which induces a serious loss of information.

The simple trace formula of Arthur ([Ar1]) applies to a particular class of functions as well ("elliptic" at two places) and to particular groups (including general linear groups). Unlike the simple trace formula of Kazhdan and Deligne, this class of functions does not automatically kill the trace of any non cuspidal representation. So it gives important information on some class of discrete series which are not cuspidal (it has been used in [Ba2], for example, for treating some residual representations).

In the next sections we will see parts of the full trace formula (for any $f\in \H(G(\a))$, any $f'\in \H(G'(\a))$) specified to general linear groups as it appears in [JL] and [AC]. They are more complicated and involve representations which are not discrete series. But even those formulae are simple compared to the full trace formula for other reductive groups (than general linear).

It is hard to give one reference for the full trace formula, because the matter is extremely difficult and occupies many long and difficult articles of James Arthur. So we send the reader to Arthur's paper [Ar2] which is the most clear possible report and contains the references. 

When using simplified trace formulae, one may get only a partial global correspondence at Step 5. However, that may be enough to accomplish a local correspondence as in Step 6 (for example, the proof of the local Jacquet-Langlands correspondence in [Ro1] and in [DKV] uses the simple trace formula of Kazhdan and Deligne and follows the described steps).\\
\ \\
\ \\
{\bf Step 3?} The comparison term by term of the geometric side of the trace formulae as given here is correct. But, as a full trace formula is more complicated that the one indicated here at Step 2, which is true only in the compact quotient case, more terms appear on the geometric side. And the other terms appearing are much more difficult to transfer between the two groups.

Also, here, we tacitly assume that if a global elliptic element $\g_o$ does not transfer (from $G(F)$ to $G'(F)$), then there is a place $w$ such that $\g_{o,w}$ does not transfer (from $G_{w}$ to $G'_w$). This is more complicated than the local transfer, and involves class field theory results -- see Lemma 4.1 from [BaRo].\\
\ \\
\ \\
{\bf Step 4?} There are three problems about claiming "By linear independence of characters on the group $G^S$...". 

First, linear independence of characters is not something obvious and might very well be untrue for some groups or for some type of representations. It is known to be true for admissible representations on reductive groups over local fields ([Cas]), but here we deal with a group which is a restricted product of a non finite number of reductive groups over different local fields.

Secondly, the number of representations appearing here is infinite, even if it might be finite for every given pair of functions $f\lra f'$.

Also, usually, a linear independence of characters might apply to an equality true for every possible good function (for example, for a reductive group over a local field, equality is required for every function in the Hecke algebra, or at least every function with support in the regular semisimple set), while our functions here have to be spherical at almost every place.\\

But it is possible to give a meaning to that linear independence of characters and to prove it. It is exposed by Flath in [Fl2] and I detailed the beautiful proof in the Appendix of the present paper in the more general setting of reductive groups. 
So this step is essentially true the way it was stated here. However, one has to fix first a representation $\Pi'$ of $G'(\a)$ and to enlarge the set $S$ such that it contains not only the places where $G'$ does not split, but also $V_\infty$, all the finite places $v$ such that $\H(G_v//K_v)$ is not commutative (see Chapter 1, Section \ref{heckealg}) and all the places where $\Pi'$ is not unramified. 

Then, let $t={\Pi'}^S$. If $f^S:=\otimes_{v\in V\bc S}f_v$, we have $\tr(t(f^S))=\prod_{v\in V\bc S}\tr\Pi'_v(f_v)$, where all the $\Pi'_v$ are unramified (by choice of $S$). We put in the equality 
$$\sum_{i\in I} m_i\tr\pi_i(f)\ =\ \sum_{j\in J} m'_j\tr\pi'_j(f')$$
only functions such that $f_v$ is spherical when $v\in V\bc S$. That kills any representation $\pi_i$ or $\pi'_j$ which is not unramified at every place $v\in V\notin S$ (Chapter 1, Lemma \ref{spherical}). So, if we stay with only this type of functions we have
$$\sum_{i\in I^0} m_i\tr\pi_i(f)\ =\ \sum_{j\in J^0} m'_j\tr\pi'_j(f')$$
where $I^0$ is the set of $i$ such that $\pi_{i,v}$ is unramified for every $v\notin S$ and  $J^0$ is the set of $j$ such that $\pi'_{j,v}$ is unramified for every $v\notin S$. Then, using the Lemma in the Appendix, we end up with 
$$\sum_{i\in I_t} m_i\tr\pi_i(f)\ =\ \sum_{j\in J_t} m'_j\tr\pi'_j(f')$$
where $I_t$ (resp. $J_t$) is the set of representations $\pi_i$, $i\in I$ (resp. $\pi'_j$, $j\in J$), such that $\pi_i^S\simeq t={\Pi'}^S$ (resp. ${\pi'_j}^S\simeq t={\Pi'}^S$). Now $\Pi'$ belongs to $J_t$. We may go and search a representation in $I_t$ corresponding to $\Pi'$.\\
\ \\
\ \\
{\bf Step 5?} The strong multiplicity one theorem is indeed true for general linear groups. It is untrue for other reductive groups, for example for $SL_n$ if $n\geq 3$ ([Bl]), and the correspondences for groups other than linear may not be injective maps, but more something like correspondences packet to packet. However, the strong multiplicity one theorem for linear groups over division algebras $D\neq F$ is not known a priori, but only {\it after} having proved a global correspondence with $GL_n(F)$ so one has to be more tricky in the proof than indicated at Step 5. The strong multiplicity one theorem for $GL_n(F)$ has been proved in [Sh] and [P-S] using Whittaker models which are not available for example on general linear groups over division algebras $D\neq F$. So the equality we get is not 
$$\tr\Pi(f)=\tr\Pi'(f')$$
if $f\lra f'$, but, as we may use the strong multiplicity one only on the left, the equation resembles more to something like 
$$\tr\Pi(f)=\sum_{j\in J_t}m'_j\tr\Pi_j'(f')$$
if $f\lra f'$. In [Ba3] it is shown that, when $G'=GL_n(D)$, the cardinality of $J_t$ is finite. This is essential for the proof of a local correspondence for square integrable representation, as in [DKV], when orthogonality of character is used. When using the general trace formula we do not have only  discrete series in the formula, but also representations induced from discrete series of proper Levi subgroups.

For other groups, one expects a sum on the left and a sum on the right and no possibility to go further (correspondence packet to packet).
\ \\
\ \\
{\bf Step 6?}
The problem with the Step 6 is that even if the Step 5 were perfect, i.e. we could come to an equality $\tr\Pi(f)=\tr\Pi'(f')$, not every unitary representation $\pi'$ may be realized as the local component of a discrete series $\Pi'$ of $G'(\a)$. So even in the best case, one cannot obtain a correspondence for unitary representations just from the trace formula. Who may be realized as a local component of a discrete series? 
It is known that the claim "any tempered representation may be realized as a local component of a discrete series" fails for any group possessing proper parabolic subgroups.
It is also known for a general reductive group that every unitary cuspidal representation may be realized as a local component of a discrete series. More generally, given local irreducible representations $\pi_{v_i}$ of $G_{v_i}$, $1\leq i\leq k$, {\it which are all unitary cuspidal}, at a finite number of places $v_1, v_2,...,v_k$, there exists a global discrete series whose local component at the place $v_i$ is $\pi_{v_i}$, $1\leq i\leq k$ ([He], Appendix 1). In the particular case of the general linear group over $F$, that is true replacing {\it unitary cuspidal} with {\it square integrable}, which is a substantial improvement. That allows one to prove the local correspondence (see [DKV]) for square integrable representations. Then the correspondence has to be extended locally to unitary representations by other means ([Ta1], [BHLS]).\\
\ \\
{\bf Back to Steps 5 and 6.} For the Jacquet-Langlands correspondence, one has to apply first the simple trace formula of Kazhdan and Deligne to get a very partial global correspondence and from that a local correspondence for square integrable representations ([DKV]). Because not every local unitary representation appears in the trace formula, one has then to extend the local correspondence to all unitary representations by local means (classification of representations and character formulas). Only then one has to come back to the global correspondence: apply the full trace formula, and use the Lemma in the Appendix to simplify with the ${}^S-part$ to get a relation of type:
$$\tr\Pi(f)=\sum_{j\in J_t}m'_j\tr\Pi'_j(f')$$
(actually more complicated, as we will see), with for example $t={\Pi'}^S$, for a discrete series $\Pi'$ of $G'$ and a well chosen finite set of places $S$. Now after simplification (i.e. fix good spherical functions at every place $v\notin S$) we end up with 
$$\prod_{v\in S} \tr \Pi_v(f_v)=\sum_{j\in J_t}m'_j\tr\Pi'_{j,v}(f'_v)$$
if $f_v=f'_v$ for $v\in S$ such that $G'_v=G_v$ and $f_v\lra f'_v$ if $G'$ does not split at $v$. Now we transfer $\Pi_v$ on $G'_v$ by the local Jacquet-Langlands correspondence, and we find an irreducible representation of $G'_v$ up to a sign. We apply the linear independence of characters on $G'_S$ and we obtain that on the right side there was only one representation, hence $\Pi'$, with multiplicity one, and such that $\Pi_v$ corresponds to $\Pi'_v$ for every $v\in S$. We applied the linear independence of characters on functions with support included in the semisimple regular set, which is not exactly the result we quoted at Section \ref{heckealg} of Chapter 1. It is a deep result of Harish-Chandra that one implies the other.\\
\ \\

The reader will find a more precise discussion based on examples in the following sections where "true" proofs are given.

\section{The case $n=1$, $d=2$}

This is the original case treated by Jacquet and Langlands in [JL], written also by Gelbart and Jacquet in [GJ]. I will follow here [GJ]. 

Let $F$ be a number field and set $G:=GL_2(F)$. Let $D$ be a central division algebra over $F$ of dimension $4$. Set $G'=D^\t$. We use the previous notation, for example $V$ is the set of places of $F$, $V_\infty$ the set of infinite places and so on. The global situation is pretty easy in this case, because, for $v\in V$, $G'_v$ is either isomorphic to $GL_2(F_v)$, or isomorphic to $D_v^\t$ where $D_v$ is a central division algebra of dimension $4$ over $F_v$. In the first case all the transfer between $G_v$ and $G'_v$ is trivial, in the second case $G'_v$ is compact modulo its center, has no proper Levi subgroup and its harmonic analysis is simple.  We follow our steps 1-5 along the proof in [GJ].\\
\ \\
{\bf Step 1.} The transfer has been explained in some detail in our paper in general (any $n$, any $d$). But the transfer of conjugacy classes (which is the basis of everything) in the case $n=1,\ d=2$ works over every local field $F_v$ in the same way as when $F_v=\r$,  
 $G_v=GL_2(\r)$ and $G'_v=\h^\t$, where $\h$ is the quaternion algebra over $\r$. We chose to explain this case because even readers which are not familiar with $p$-adic fields may easily follow. Here $(1,i,j,k)$ is the usual basis of $\h$ and $a,b,c,d\in\r$. If $x=a+bi+cj+dk\in \h$, define the trace and the norm of $x$ by  $\tr(x)=2a$ and $N(x)=a^2+b^2+c^2+d^2$ and set $P_x(X)=X^2-\tr(x)X+N(x)$. Using the formula $N(x)=(a+bi+cj+dk)(a-bi-cj-dk)$ it is easy to check that $P_x(x)=0$. 
If $x\notin Z=\r$, then $\tr(x)^2<4N(x)$ and $P_x$ is irreducible and is the minimal polynomial of $x$ over $\r$. So the regular semisimple elements of $\h^\t$ are the elements in $\h^\t\bc Z^\t$.
Moreover, if $Q(X)=X^2-uX+v$ is an irreducible element of $\r[X]$, then $u^2<4v$ so there is $x\in{\h^\t}^{rs}$ such that $Q$ is the characteristic (and minimal) polynomial of $x$. Then $x$ corresponds to $g$, where  $g\in GL_2(\r)$ is
$$
g=\begin{pmatrix}
0&-v\\
1&u
\end{pmatrix}.$$ 
The elements of $GL_2(\r)$ which are not in the image of the transfer are the matrices with real eigenvalues. The situation is the same replacing $\r$ by $F_v$. {\it However, the parallel with the situation over $\r$ applies only in our case $n=1$, $d=2$, because there is no central division algebra of dimension superior to $4$ over $\r$, whereas such algebras always exist when $F_v$ is a $p$-adic field.}\\
\ \\
{\bf Step 2.} The trace formulae for the two groups are computed in [GJ]. We explain the link with our computation for the compact quotient case in Section \ref{compactquotient}. We adopt, as in [GJ], the notation $\tilde{G}$ and $\tilde{G'}$ for $Z\bc G$ and $Z\bc G'$ where $Z$ is the center.

Because $\tilde{G}'(F)\bc \tilde{G}'(\a)$ is compact, the trace formula for $G'$ is the one we computed in Section \ref{compactquotient}:

$$\tr(R'(f'))=\sum_{\oo'\in X'} \vol(\tilde{G}'_{\g_{o'}}(F)\bc \tilde{G}'_{\g_{o'}}(\a))\Phi(f'_\o,\g_{o'}).$$

This is the formula of (1.11) in [GJ]. Theorem 1.16 [GJ] claims that $R'(f')$ is trace class (their reference is [DL]). The type of functions they use are functions with compact support modulo the center on which the center acts by  the character $\o^{-1}$, and this corresponds to our functions $f_\o$, $f'_\o$, not $f$, $f'$. So our $f'_\o$ here is the one denoted $\phi$ in [GJ]. Notice that $\pi(f_\o)$ is defined only for representations $\pi$ with central character $\o$, by the formula $\pi(f_\o):=\int_{Z(\a)\bc G(\a)}f_\o(\bar{g})\pi(\bar{g})\ d\bar{g}$, and so we have $\pi(f_\o)=\pi(f)$.\\

The trace formula for $G$ is much more complicated. I give here the result of Theorem 6.33 in [GJ]. The representation  $R$ in the space $L^2(G(F)Z(\a)\bc G({\a}))$ as denoted in our paper for the compact quotient case is no longer the good one as it contains a continuous part. The trace formula involves its restriction to the space $L^2_{disc}(G(F)Z(\a)\bc G({\a}))=L^2_0(G(F)Z(\a)\bc G({\a}))\oplus L^2_{res}(G(F)Z(\a)\bc G({\a}))$, i.e. to the sum of irreducible subrepresentations of the space $L^2(G(F)Z(\a)\bc G({\a}))$. It is this representation I will denote here by $R$. In [GJ] it is denoted $\rho_{cusp}+\rho_{sp}$. Also, $X$ will not denote here the set of all conjugacy classes of $\tilde{G}(F)$, but only the set of {\it elliptic} classes, i.e. the trivial class and those classes for which the characteristic polynomial of a representative in $G(F)$ is irreducible 
(this is well defined as the property is $Z$-invariant). The function $\phi$ in [GJ] is our function $f_\o$ here. Corollary 2.4 [GJ] claims that $R(f)$ is trace class (there is a proof, using again [DL]), and the formula ([GJ] page 240) is:

$$\tr(R(f))=\sum_{\oo\in X} \vol(\tilde{G}_{\g_{o}}(F)\bc \tilde{G}_{\g_{o}}(\a))\Phi(f_\o,\g_{o})$$
$$ (6.34)\ \ \ \ \ \ \ \ \ \ \ \ \ \ \ \ \ \ \ \ \ \ \ \ \ \ \ \  \ \   +\ {\rm f.p.} \int _{Z(\a)N(\a)\bc G(\a)}  f_\o(g^{-1}
\begin{pmatrix}
1&1\\
0&1
\end{pmatrix}
g)\ dg$$
$$(6.35)\ \ \ \ \ \  -\ \frac{1}{2}vol(F^*\bc F^0(\a))\int_K\int_{N(\a)}\sum_{\alpha\neq 1} f_\o(k^{-1}n^{-1}
\begin{pmatrix}
\alpha & 0\\
0&1
\end{pmatrix}
nk){\rm log} H[wnk]\ dn\ dk$$
$$ (6.36)\ \ \ \ \ \ \ \ \ \ \ \ \ \ \ \ \ \ \ \ \ \ \ \ \ \ \ \  \ \ \ \ \ \   +\ \frac{1}{4\pi}
\int_{-\infty}^{\infty} \tr(M(-iy)M'(iy)\pi_{iy}(f_\o))\ dy$$
$$ (6.37)\ \ \ \ \ \ \ \ \ \ \ \ \ \ \ \ \ \ \ \ \ \ \ \ \ \ \ \  \ \ \ \ \ \ \ \ \ \ \ \ \ \ \ \ \ \ \ \ \ \ \ \ \ \  -\ \frac{1}{4} \tr(M(0)\pi_{0}(f_\o)).$$

We do not explain the terms yet.\\
\ \\
{\bf Step 3.}
Now Gelbart and Jacquet work out the formula of $R(f)$ in order to show that when $f\lra f'$ some terms are zero and to compare it with the formula for $R'(f')$.\\ 

The term (6.34):  one has
$$ \int _{Z(\a)N(\a)\bc G(\a)}  f_\o(g^{-1}
\begin{pmatrix}
1&1\\
0&1
\end{pmatrix}
g)\ dg=
\prod_{v\in V}\int _{Z_vN_v\bc G_v}  {f_{\o,v}}(g_v^{-1}
\begin{pmatrix}
1&1\\
0&1
\end{pmatrix}
g_v)\ dg_v.
$$
We assume (without loss of generality, see Section \ref{transfct}) that the support of $f_v$ is included in $G_v^{rsd_v}$ if $G'$ does not split at $v$ (i.e. $d_v=2$). (This assumption is not made in [GJ].)
Then the elements conjugated to 
$\begin{pmatrix}
1&1\\
0&1
\end{pmatrix}$ 
are not in the support of ${f_{\o,v}}$ (because their characteristic polynomial is not irreducible) and so the term (6.34) is zero. In [GJ] the assumption on the support is not made, that is why their computation is more complicated and involves the $L$-functions.
Our assumption has a price to pay, but in characteristic zero it is very little, since it is proved that the harmonic analysis of $G_v$ and $G'_v$ are ruled by functions with regular support.\\ 

The term (6.35): same proof : if $v$ is a place where $G'_v$ does not split, then elements conjugated to 
$\begin{pmatrix}
\alpha & 0\\
0&1
\end{pmatrix}$
are not in the support of $f_{\o,v}$ because they have real eigenvalues (see Step 1).\\

The term (6.36): it is shown in [GJ] page 243 that this term equals (see their formula (7.13)):
$$(7.13)\ \ \ \ \frac{1}{4\pi}\sum_{\chi}\int_{-i\infty}^{i\infty} m'(\eta)\tr(\pi_{\eta}(f_\o))$$
$$+\frac{1}{4\pi} \sum_{u\in V}\sum_{\chi}\int_{-i\infty}^{i\infty} \prod_{v\in V\bc u} \tr(\pi_{\eta_v}(f_{\o,v}))\tr(R_u(\eta_u)^{-1}R'(\eta_u)\pi_{\eta_u}(f_{\o,u}))\ dy.$$
Without getting into details, let us say that $\pi_\eta$ is a representation of $G(\a)$ parabolically induced from a representation of the diagonal torus which depends on the character $\eta$, which depends itself on the character $\chi$ and on the complex number $y$. The only thing that matters for us is that, if $v\in V$ is such that $G'_v$ does not split, then $\tr(\pi_{\eta,v}(f_{\o,v}))=0$ by Proposition \ref{inductionjl} (a) since the diagonal torus of $G_v$ does not correspond to a Levi subgroup of $G'_v$, (because it is compact modulo its center, $G'_v$ has no proper Levi subgroup). As there is at least one place $v$ where $G'_v$ does not split, the first term of (7.13) is zero. Now the second term of (7.13): we will show that every product appearing is zero. It is enough to show that there are at least two places $v_1,v_2\in V$ such that $G_{v_1}'$ and $G'_{v_2}$ do not split, because then $\tr(\pi_{v_1}(f_{\o,v_1})) = \tr(\pi_{v_2}(f_{\o,v_2})) = 0$ so each product contains a zero term. The fact that two such places exist, i.e. the number of places where $D$ does not split is at least two, is a classical  (non trivial) result from the theory of central simple algebras over number fields, which is part of the class field theory. {\it This term resembles to the splitting formulae for $(G,M)$-families, which appear in the work of Arthur for general groups 
{\rm (}see for example {\rm Lemma 17.6} of {\rm [Ar2]} -- there is a typo in formula {\rm (17.14)}, the last $Q_1$ is a $Q_2${\rm )}. Sometimes, after the splitting formula has been applied, the same type of game between two different places combined with class field theory is played to show that the term vanishes {\rm (}{\rm Corollary 7.5} in {\rm [Ar1]} or {\rm Lemma 4.5} and following in {\rm [BaRo]} for example use such proofs{\rm )}}.\\

The term (6.37): the operator $M(0)$ intertwines the representation $\pi_0$ with itself ([GJ] (3.22) for $s=0$). But $\pi_0$ is irreducible because it is induced from a unitary irreducible representation of the diagonal torus. So $M(0)$ is a scalar. Then, in order to show that the term (6.37) is zero, it is enough to show that $\tr(\pi_0(f_\o))=0$. This is true by Proposition \ref{inductionjl} (a), because there exists a place $v\in V$ such that $G'_v$ does not split, and at this place the diagonal torus is a Levi subgroup which does not transfer. {\it In the general case $n\geq 2$, after comparison of trace formulae, {\rm [AC]} got rid of all the other terms but those ones, which are not zero like here; they appear in the following section, where the reader will encounter some coefficients $\frac{1}{l^2}$, where $l$ is the number of blocks of a Levi subgroup, of which the $\frac{1}{4}$ here is the particular case $l=2$}.\\
\ \\

Now that we showed that if $f\lra f'$ we have:
$$\tr(R(f))=\sum_{\oo\in X} \vol(\tilde{G}_{\g_{o}}(F)\bc \tilde{G}_{\g_{o}}(\a))\Phi(f_\o,\g_{o}),$$
and
$$\tr(R'(f'))=\sum_{\oo'\in X'} \vol(\tilde{G}'_{\g_{o'}}(F)\bc \tilde{G}'_{\g_{o'}}(\a))\Phi(f'_\o,\g_{o'}),$$
we want to show that the right parts of the formulae are equal.

The terms corresponding to $\g=1$ and $\g'=1$ are zero since we have functions which are supported in the regular semisimple set at least at one place, and $1$ is not regular semisimple. If $\g\neq 1$, $\g'\neq 1$ and $\g\lra \g'$, then the terms corresponding to $\g$ and $\g'$ are equal by definition of $f_\o\lra f'_\o$ and the choice of measures.
It is enough then to show that the terms indexed by elements of $X\bc \{1\}$ which do not correspond to elements of $X'$ are zero. This is a consequence of the class field theory as explained in the next paragraph.

Let $\g$ be an element of $GL_2(F)$ such that the characteristic polynomial $P_\g$ of $\g$ is irreducible over $F$ (i.e. the class of $\g$ is in $X\bc\{1\}$). Assume $\Phi(f_{\o},\g)\neq 0$. 
Because  
$\Phi(f_{\o},\g)=\prod_{v\in V} \Phi(f_{\o,v},\g_{v})$, 
we have that for every place $v$ such that $G'_v$ does not split $\g_{v}$ corresponds to an element of $G'_v$ i.e. $P_\g$ (which is the characteristic polynomial of $\g_{v}$) does not split over $F_v$. Let $L$ be the subfield of $M_2(F)$ generated by $\g$, $L\simeq F[X]/(P_\g)$. Then, for every place $v$ such that $G'_v$ does not split, $L_v:=L\otimes F_v$ is a {\it field} extension of degree $2$ of $F_v$, because it is isomorphic to $F_v[X]/(P_{\g})$. The proposition 5 of [We], XIII.3, more precisely "(ii) implies (iii)", shows then that there is an $F_v$-algebras isomorphism from $L$ to a subfield of $D$. Then the image $\g'$ of $\g$ is an element of $D$ having the same characteristic polynomial as $\g$. So, if $\Phi(f_{\o},\g)\neq 0$, then there is $\g'\in G'(F)=D^\t$ such that $\g\lra \g'$.\\
\ \\
{\bf Step 4.} We start with the equality $\tr(R(f))=\tr(R'(f'))$.  
The Corollary 1.7 and Corollary 2.5 in [GJ] claim that $R$ and $R'$ decompose in a sum of irreducible representations with finite multiplicities (the references are [GGP-S] and [La]). We write 
$$\sum_{i\in I}m_i\tr\pi_i(f_\o)=\sum_{j\in J}m'_j\tr\pi'_j(f'_\o)$$
for $f\lra f'$.
Recall $\tr\pi(f_\o)=\tr\pi(f)$ and $\tr\pi'_j(f'_\o)=\tr\pi'_j(f')$  so we simply write:
$$\sum_{i\in I} m_i\tr\pi_i(f) = \sum_{j\in J} m'_j\tr\pi'_j(f')$$
for $f\lra f'$. {\bf{We fix a representation $\pi'_{j_0}$ on the right}}. Let $S$ be the set of places $v$ of $F$ such that either $v$ is infinite, or $G'$ is not split at $v$, or $\pi'_{j_0,v}$ is not unramified. 

Now every $\pi_{i,v}$ and every $\pi'_{j,v}$ is unitary. For each collection  $u:=\{(u_v)\}_{v\in V\bc S}$ where, for every $v\in V\bc S$, $u_v$ is an irreducible unitary class of representations of $G_v=G'_v$, denote by $A(u)$ the subset of $I$ made of elements $i$ such that $\pi_{i,v}\in u_v$ for all $v\in V\bc S$,  and by $A'(u)$ the subset of $J$ made of elements $j$ such that $\pi'_{j,v}\in u_v$ for all $v\in V\bc S$. Let $U(S)$ be the set of $u$ such that $A(u)$ or $A'(u)$ is not empty. 

If $f\lra f'$ such that $f_v$ is a spherical function on $G_v=G'_v$ for every $v\in V\bc S$ we have 

$$\sum_{u\in U(S)} c_u \prod_{v\in V\bc S} \tr u_v(f_v)=0,$$
where 
$$c_u= \sum_{i\in A(u)}\prod_{v\in S}\tr\pi_{i,v}(f_v)\ -\ \sum_{j\in A'(u)}\prod_{v\in S}\tr\pi'_{j,v}(f'_v).$$ 
For every $u$ such that there is a $v\in V\bc S$ such that $u_v$ is not unramified, we have $\tr u_v(f_v)=0$ because $f_v$ is spherical (the trace of a representation $\pi$ on a spherical function is zero unless $\pi$ is unramified, Chapter 1, Lemma \ref{spherical}). After deleting those terms we remain only with elements $u$ of $U(S)$ such that $u_v$ is unramified for every $v\in V\bc S$. Applying the result from the Appendix, we conclude that $c_u=0$ for those $u$ ([GJ] has an alternate approach here, based on [LL]; we chose Flath's proof which is reproduced in our Appendix).
But one of those $u$ is $\{\pi_{j_0,v}\}_{v\in V\bc S}$. We find 
$$\sum_{i\in A(\pi'_{j_0})}\ m_i\prod_{v\in S}\tr\pi_{i,v}(f_v) = \sum_{j\in A'(\pi'_{j_0})}\ m'_j\prod_{v\in S}\tr\pi'_{j,v}(f'_v)$$
where $A(\pi'_{j_0})$ is the set of $i$ such that $\pi_{i,v}\simeq \pi'_{j_0,v}$ for all $v\in V\bc S$ and 
$A'(\pi'_{j_0})$ is the set of $j$ such that $\pi'_{j,v}\simeq \pi'_{j_0,v}$ for all $v\in V\bc S$. Here $f'_v=f_v\in \H(G'_v)=\H(G_v)$ for the places $v\in S$ such that $G'_v=G_v$ and $f_v\lra f'_v$ for the places $v$ where $G'_v$ does not split.\\ 
\ \\
{\bf Step 5.} Now we have the following {\it finitude theorems}:\\

- $A(\pi'_{j_0})$ has at most one element, which appears with multiplicity $1$ by the strong multiplicity one theorem. ([P-S], [Sh])
\\

- $A'(\pi'_{j_0})$ is finite by Lemma 5.14 in [Ro1].\\

We have moreover linear independence of characters on the groups $G'_v$, $v\in S$. When $G'_v$ does not split, $G'_v$ is compact modulo its center. It is easy to show then that if a linear combination of traces of irreducible representations is zero on every $f'\in \H(G'_v)$ supported in the regular set, then it is zero on every function in $\H(G'_v)$\footnote{When $n>2$, $G'_v$ is no longer compact modulo its center, but the same holds; this is a deep result of Harish-Chandra.}. So, because on the right hand side of the equality at least the one representation $\pi'_{j_0}$ appears, the sum is not identically zero, which shows that $A(\pi'_{j_0})\neq \emptyset$. Then  $A(\pi'_{j_0})$ has one element, say $\pi_0$, which already has the property that it has local components $\pi_{0,v}\simeq \pi'_{j_0,v}$ for all places $v\in V\bc S$. Moreover, we have the equality
$$\prod_{v\in S}\tr\pi_{0,v}(f_v) = \sum_{j\in A'(\pi'_{j_0})}\ m'_j\prod_{v\in S}\tr\pi'_{j,v}(f'_v),$$
where $A'(\pi'_{j_0})$ is finite. Now $S$ contains the set $S_0$ of places $v$ where $G'_v$ does not split. By linear independence of characters on the groups $G'_v=G_v$ for $v\in S\bc S_0$, for every $j\in A'(\pi'_{j_0})$ we have $\pi'_{j,v}\simeq \pi_{0,v}$ for all $v\in S\bc S_0$, and simplify:
$$\prod_{v\in S_0}\tr\pi_{0,v}(f_v) = \sum_{j\in A'(\pi'_{j_0})}\ m'_j\prod_{v\in S_0}\tr\pi'_{j,v}(f'_v).$$

In order to prove the global correspondence now, {\it we will show that on the right-hand side the only representation is $\pi'_{j_0}$} and so $\pi_{0,v}$ corresponds to $\pi'_{j_0,v}$ for the places $v\in S_0$, where $G'_v$ does not split (we know already that they are isomorphic at all the places not in $S_0$).

The classification of discrete series of $GL_2(\a)$ as proved in [MW1], implies that {\it $\pi_0$ is either a unitary character, or a cuspidal representation}.\\
\ \\
{\bf First case: $\pi_0$ is a unitary character.} In this case we will directly construct a character $\chi'$ of $G'(\a)$ corresponding to $\pi_0$. It is known (classical result for $GL_n$) that a smooth unitary character of $GL_2(\a)$ is of the form 
$$\chi:g\mapsto |\det(g)|_\a^{i\alpha_\chi}=\prod_{v\in V} |\det(g_v)|_v^{i\alpha_\chi},$$ 
where $g=(g_v)_{v\in V}$, $\alpha_\chi$ is a real number and $|\ |_v$ is the normalized norm of the local field $F_v$. 
 
If $v\in S_0$ and $G'_v=D_v^\t$, where $D_v$ is a quaternion division algebra over $F_v$, a multiplicative norm $N_v:G'_v\to F_v^\t$ may be defined as in the case $F_v=\r$ treated at Step 1. 
If $g'_v\in {G'}_v^{rs}$, then $N_v(g'_v)$ is the constant coefficient of the minimal polynomial of $g'_v$ (as in the case $F_v=\r$). So $\det(g_v)=N_v(g'_v)$ if $g_v\lra g'_v$. Then 
$$g'\mapsto \prod_{v\in S_0} |N_v(g'_v)|_v^{i\alpha_\chi}\prod_{v\in V\bc S_0} |\det(g'_v)|_v^{i\alpha_\chi}$$
is a character of $G'(\a)$ which corresponds to $\chi$. Specializing $\chi$ to $\pi_0$, there is a character $\chi'$ of $G'(\a)$ such that $\pi_{0,v}$ corresponds to $\chi_v'$ for all $v\in S_0$ and $\pi_{0,v}=\chi'_v$ for every $v\in V\bc S_0$. This second property of $\chi'$ shows that it appears on the right hand side of our equality. Because it corresponds to $\pi_0$, we may simplify the equality by dropping $\chi$ from the left and $\chi'$ from the right. But then, by linear independence of characters on the groups $G'v$, there is nothing left on the right hand side. So $\chi'=\pi'_{j_0}$.\\
\ \\
{\bf Second case: $\pi_0$ is cuspidal.} Now local information comes into play to be combined with global information. Let $v\in S_0$. Because the equality is not identically $0$, the trace of $\pi_{0,v}$ is not zero on all the functions $f$ which correspond to a function $f'$. So $\pi_{0,v}$ is not an induced representation from the unique standard Levi subgroup which is the diagonal torus of $G_v=GL_2(F_v)$ (Proposition \ref{inductionjl} (a)). The classification of irreducible representation of $G_v$ by parabolic induction implies that $\pi_{v,0}$ is either of finite dimension (actually of dimension one if $v$ is finite), or a {\it square integrable representation} (I will not define it here). Moreover, one knows that, because $\pi_0$ is cuspidal, $\pi_{0,v}$ is a generic representation for every $v\in V$. I do not define {\it generic} here, but a finite dimensional representation is not generic so the only possibility is $\pi_{0,v}$ is square integrable.\\

Now the {\it orthogonality of characters} is needed. If $\pi'$ is an admissible irreducible representation of $G'_v$, then $\pi'$ has a central character $\o_v$. Because $G'_v$ is compact modulo the center, $\pi'$ is finite dimensional and one may define the {\bf function character}  $\chi_{\pi'}$ of $\pi'$ by $\chi_{\pi'}:G'_v\to \cc$, $\chi_{\pi'}(g')=\tr\pi'(g')$. Then $\chi_{\pi'}$ is locally constant, it is a class function (i.e. stable by conjugation) and has the same central character as $\pi'$ (i.e. $\chi_{\pi'}(zg')=\o_v(z)\chi_{\pi'}(g')$ for all $z\in Z_v$ and all $g'\in G'_v$). The link between the distribution character and the function character is : $\tr{\pi'}(f')=\int_{G'_v}\chi_{\pi'}(g')f'(g')dg'$ if $f'\in\H(G'_v)$. 

If $T'$ is a maximal torus of $G'_v$, then $T'/Z_v$ is compact, and on $T'$ we chose the Haar measure such that the volume of $T'/Z_v$ is one. This verifies that when two tori are conjugate, the measures are conjugate. Let ${\bf S}$ be a set of representatives of conjugacy classes of maximal tori of $G'_v$. For each $T'\in {\bf S}$ let $|W_{T'}|$ be the cardinality of the Weyl group of $T'$. Then there is a scalar product defined on the set $C(G'_v,\o_v)$ of continuous class functions on $G'_v$ with central character $\o_v$ by the formula:
$$(f,h):=\sum_{T'\in {\bf S}}\ \frac{1}{|W_{T'}|}\ \int_{T'/Z_v} \overline{f({t})} h({t}) D({t})dt$$
where $D$ is a class function with values in $\r_+^*$ which depends only on the eigenvalues of $t$. 

We extend the product to the space $L^2(G'_v,\o_v)$ of class functions $f$ with central character $\o_v$ such that $(f,f)$ converges. Then the set $\{\chi_{\pi'}\}$ where $\pi'$ runs over the set of isomorphism classes of admissible irreducible representations of $G'_v$ with central character $\o_v$ form a orthonormal subset of $L^2(G'_v,\o_v)$. This is a classical result in the theory or representations of compact groups, where the scalar product is taken as simply the integral on the group, combined with the Weyl integration formula, which takes a particular form when the functions we integrate are class functions. This form of the scalar product is useful because we do not have point by point transfer of functions from $G_v$ to $G'_v$, but we have class functions transfer, as we transferred the conjugacy classes.

Let now $\pi$ be a smooth irreducible representation of $G_v$. According to [H-C2], there is a unique locally constant function $\chi_\pi:G_v^{rsd_v}\to \cc$ (recall here $d_v=2$), such that:

- $\chi_\pi$ is stable by conjugation with elements of $G_v$,

- $\chi_\pi$ has central character $\o_v$ (the set $G^{rsd_v}$ is stable by $Z_v$),

- $\tr\pi(f)=\int_{G_v^{rsd_v}}\chi_\pi(g) f(g) dg$ for all $f\in \H(G_v)$ supported in $G_v^{rsd_v}$.\\

Using our transfer of regular semisimple conjugacy classes (Step 1), we see every class function on $G_v^{rsd_v}$ as a class function on ${G'}_v^{rs}$. 
Deep results are that, if $\pi$ is square integrable, then $\chi_\pi\in L^2(G'_v,\o_v)$ (notice that the set $G_v'\ \bc\ {G'}_v^{rs}$ is of measure zero) and the norm of $\chi_\pi$ is $1$ [Cl2]. Let us now get back to the equation:

$$\prod_{v\in S_0}\chi_{\pi_{0,v}}(g_v) = \sum_{j\in A'(\pi'_{j_0})}\ m'_j\prod_{v\in S_0}\chi_{\pi'_{j,v}}(g'_v)$$
if $g_v\lra g'_v$ for $v\in S_0$. Computing the square of the norm of the left and right hand sides (one uses the isomorphism between $L^2(X\t Y)$ and $L^2(X)\otimes L^2(Y)$ for $X$ and $Y$ measured spaces, since the scalar product has been defined so far only place by place) we find
$$\prod_{v\in S_0}||\chi_{\pi_{0,v}}||^2 = \sum_{j\in A'(\pi'_{j_0})}\ {m'}^2_j\prod_{v\in S_0}||\chi_{\pi'_{j,v}}||^2.$$
All the norms appearing are $1$ (we are in the case where $\pi_{0,v}$ are square integrable for $v\in S_0$); we have then:
$$1=\sum_{j\in A'(\pi'_{j_0})} {m'}_j^2.$$
Because $m'_j$ are integers, on the right there was only one representation, $\pi'_{j_0}$, and its multiplicity is $1$. Now the relation
$$\prod_{v\in S_0}\tr\pi_{0,v}(f_v) = \prod_{v\in S_0}\tr\pi'_{j_0,v}(f'_v)$$
implies that for every $v\in S_0$, there is a scalar $\l_v\in \cc^\t$ such that $\pi_{0,v}$ corresponds to $\l_v \pi'_{j_0,v}$, so $\chi_{\pi_{0,v}}=\l_v\chi_{\pi'_{j_0,v}}$. Because the norm of $\chi_{\pi_{0,v}}$ and $\chi_{\pi'_{j_0,v}}$ is one, $|\l_v|=1$. To get to our definition of the global correspondence, we have moreover to prove that $\l_v=\pm 1$, or, equivalently, $\l_v\in\r$. This is not obvious. There are two (not simple) methods. Fix some place $v_0\in S_0$. In order to show that $\l_{v_0}\in\r$, he first method is to show that there exists $g\lra g'$ such that $\chi_{\pi'_{j_0,v_0}}(g')$ and $\chi_{\pi_{0,v_0}}(g)$ are both non zero real numbers. This method is used in [Ro1], Proposition 5.9. Another method is the following: define another discrete series $\pi'$ of $G'$
such that $\pi'_{v_0} = \pi'_{j_0,v_0}$ and, for every place $v\in S_0\bc\{v_0\}$, $\pi'_{v}$ is the trivial representation of $G'_v$ (this is possible since the trivial representation and $\pi_{j_0,v_0}$ are square integrable representations; see {\bf Step 6?}). Assume we did again the whole work with $\pi'$ in place of $\pi'_{j_0}$. As the trivial representation of $G'_v$ for $v\in S_0$ corresponds to the Steinberg representation of $G_v$, we know that at these places $\l_v=-1$. Because $\prod_{v\in S_0}\l_v=1$, $\l_{v_0}=\pm 1$. 

We define the global Jacquet-Langlands correspondence by setting $\pi'_{j_0}\mapsto \pi_0$. The fact that it is injective is obvious: if $\pi'$ and $\pi''$ have the same image, then $\tr\pi'_v = \pm \tr\pi''_v$, so $\pi'_v\simeq \pi''_v$ by linear independence of the characters.
One may also check that the image of this map is the set of discrete series of $G(\a)$ which do not transfer to zero at any place $v\in S_0$. The proof is the same, replacing at Step 4 "{\bf we fix a representation $\pi'_{j_0}$ on the right}" with "{\bf we fix a representation $\pi_0$ on the left}" and changing everything accordingly.\\
\ \\

\section{The proof of the global Jacquet-Langlands correspondence}

Let $F$ be a number field, $D$ a central division algebra of dimension $d^2$ ($d\in\n^*$) over $F$ and $S_0$ the finite set of places of $F$ where $D$ does not split, {\it which we suppose to be disjoint with the set $S_\infty$ of infinite places of $F$}. For $k\in \n^*$, let $G_k:=GL_k(F)$ and $G'_k:=GL_k(D)$. Let $DS_k$ be the set of discrete series of $G_k(\a)$, $DS'_k$ the set of discrete series of $G'_k(\a)$. 
Every group $G'_{k,v}$ is of type $GL_{kr_v}(D_v)$ where $D_v$ is a central division algebra of dimension $d_v^2$, where $d_v=\frac{d}{r_v}$ over $F_v$. $d_v=1$ if and only if $v\notin S_0$.\\
 
An ingredient for the trace formula (in the non compact quotient case) are the representations induced from discrete series. We define them here. A standard Levi subgroup $L$ of $G_k$ is given by a partition $(n_1,n_2,...n_l)$ of $k$ and we use the notation $L_v$ for the standard Levi subgroup of $G_{k,v}$ corresponding to the same partition. A standard Levi subgroup $L'$ of $G'_{k,v}$ is given by a partition $(n_1,n_2,...n_l)$ of $k$ and we use the notation $L'_v$ for the standard Levi subgroup of $G'_v$ corresponding to the partition $(r_vn_1,r_vn_2,...,r_vn_l)$ of $r_vk$. If $L$ is a standard Levi subgroup of $G$ corresponding to a partition $(n_1,n_2,...n_l)$ of $k$, then a {\bf discrete series of} $L$ is a representation 
$\rho_1\otimes \rho_2\otimes ...\otimes \rho_l$, where $\rho_i$ is a discrete series of $GL_{n_i}(\a)$. We define
$\pi:= \ind_L^G \ \rho_1\otimes \rho_2\otimes ...\otimes \rho_l$ by local induction at every place, namely $\pi_v:=\ind_{L_v}^{G_{k,v}}\ \rho_{1,v}\otimes \rho_{2,v}\otimes ...\otimes \rho_{l,v}$, which is an irreducible representation, because local induction for unitary representations preserves irreducibility (Theorem \ref{unitirred}). The same holds for $G'$.

We want to prove the global Jacquet-Langlands correspondence :

\begin{theo}
{\rm (a)} There exists an injective map ${\bf G}:DS'_k\to DS_{kd}$, such that, for all $\pi'\in DS'_k$, if $\pi:={\bf G}(\pi')$, then $\pi_v\simeq \pi'_v$ for all $v\in V\bc S$ and $\pi_v$ corresponds to $\pi'_v$ for all $v\in S$.

{\rm (b)} If $\pi\in DS_{kd}$, then $\pi$ is not in the image of ${\bf G}$ if and only if there exists a place $v\in S$ such that $\pi_v$ corresponds to zero.

{\rm (c)} The strong multiplicity one theorem is true for $G'_k(\a)$: if two discrete series have isomorphic local components at almost every place then they are equal.
\end{theo}
\ \\

Here is a discussion of the Steps 1-5 for the proof of the theorem. 
{\it We prove {\rm (a)}, {\rm (b)} and {\rm (c)} altogether by induction on $k$, so we assume in the sequel the theorem is true for every $k<n$ and we prove it for $n$. We denote all the correspondences by {\bf G} without mentioning $k$}.\\

Step 1 has already been explained in Sections \ref{transconjclass}, \ref{transtori}, \ref{transfct} and we have just stated the global correspondence above. Steps 2 and 3 are very complicated and have been carried out by Arthur and Clozel in [AC], based on previous work of Arthur. I will explain the proof starting from the equality of the spectral side of the trace formula of $G_{nd}(\a)$ and $G'_n(\a)$. This is the relation at the page 203 of [AC]:

$$\sum_L |W_0^L||W_0^{G_{nd}}|^{-1}\sum_{s\in W(\frak{a}_L^G)_{reg}} |\det(s-1)_{\frak{a}_L^{G_{nd}}}|^{-1}\tr(M(s,0)\rho_{Q,t,\mu}(f))=$$

$$=\sum_{L'} |W_0^{L'}||W_0^{G'_n}|^{-1}\sum_{s\in W(\frak{a}_{L'}^{G'})_{reg}} |\det(s-1)_{\frak{a}_{L'}^{G'_n}}|^{-1}\tr(M(s,0)\rho'_{Q',t,\mu}(f')),$$
if $f\lra f'$.

Before explaining the symbols, notice that we changed $\sum_M$ (in [AC]) to $\sum_L$ (here) because another letter $M$ appears in the formula.
Also, the $l$ and the $\sigma$ disappeared from the second term of the equality, because the formula from [AC] is more general (is written in the base change setting) and applies here with $l=1$ and $\sigma$ trivial, which may be deleted.

In this equality, $L$ runs over the set of Levi subgroups of $G_{nd}$ which contain the diagonal torus and verify some symmetries ($W(\frak{a}_L^G)_{reg}\neq \emptyset$). As explained in [AC], the contribution of such a Levi subgroup to the formula is the same as the contribution of a standard Levi subgroup to which it is conjugate. The symmetries imply that the standard Levi subgroup is made of blocs of the same size. 
$W_0^L$ is the Weyl group of $L$, and $W_0^{G_{nd}}$ is the Weyl group of $G_{nd}$. 
For a standard Levi subgroup with $l$ blocs ($l$ divides $nd$), $s$ runs through the cycles of maximal length in ${\frak{S}_l}$. The contribution of any such element to the trace formula is the same as the one of the cycle $(1,2,...,l)$.
In order to not complicate things, we forget $t$ and $\mu$ and we say that $\rho_{Q,t,\mu}$ is the representation induced from the whole space of discrete series of $L$, with respect to the parabolic subgroup $Q$ of $G_{nd}$ (we have already seen that, because the induced representation is irreducible, that does not depend on $Q$, but only on $L$). $M(s,0)$ is the standard intertwining operator with respect to $s$ at $0$. The representation $\rho_{Q,t,\mu}$ is a direct sum of all irreducible representations of the form 
$\pi:= \ind_L^{G_{nd}} \ \rho_1\otimes \rho_2\otimes ...\otimes \rho_l$, where the $\rho_i$, $1\leq i\leq l$, are discrete series of the blocks of $L$. For a place $v$ where all the $\rho_i$ are unramified, $\pi_v$ is unramified (so $\pi(f)$ makes sense for $f$ in the global Hecke algebra). Now $M(s,0)$ intertwines the space of $\pi$ with the space of  $\ind_L^{G_{nd}} \ \rho_{s(1)}\otimes \rho_{s(2)}\otimes ...\otimes \rho_{s(l)}$ (as subspaces of $\rho_{Q,t,\mu}$). So, if $s$ is the cycle $(1,2,...,l)$, the contribution to the trace of the operator $M(s,0)\rho_{Q,t,\mu}(f)$ is only that of representations $\pi$ for which $\rho_1\simeq \rho_2\simeq ...\simeq \rho_l$.

The same is true for the $G'$ part of the formula. Now the representations on the $G$ side appear only once (strong multiplicity one theorem for $G_n$ and for its Levi subgroups). On the $G'$ side we have the strong multiplicity one theorem for the proper Levi subgroups by the induction hypothesis. So, after some computation ([Ba1], page 414 and the following) the equality has the form:
$$\sum_{\rho\in DS_{nd}}\tr\rho(f)\ +\ \sum_{l|nd,\ l\neq 1}\frac{1}{l^2}\sum_{\rho\in DS_{\frac{nd}{l}}} \tr(M(s_l,0)\ind(\underbrace{\rho\otimes\rho\otimes ...\otimes\rho}_{l\ {\rm times}})(f))$$
$$=\ \sum_{\rho'\in DS'_{n}} m_{\rho'} \tr\rho'(f)\ +\ \sum_{l|n,\ l\neq 1}\frac{1}{l^2}\sum_{\rho'\in DS'_{\frac{n}{l}}} \tr(M(s_l,0)\ind(\underbrace{\rho'\otimes\rho'\otimes ...\otimes\rho'}_{l\ {\rm times}})(f'))$$
if $f\lra f'$, where $s_l=(1,2,...,l)$, $m_{\rho'}$ is the multiplicity of $\rho'$ in the discrete spectrum of $G'_{n}(\a)$ and $\ind$ is the parabolic induction from the Levi subgroup $L_l$ with $l$ blocks of equal size. $M(s,0)$ intertwines the induced representation with itself and, as the induced representation is irreducible, $M(s,0)$ acts as a scalar. Because $M(s,0)$ is a unitary operator, the scalar has modulus one.

If $l$ divides $nd$ but $l$ does not divide $n$, then there is a place $v$ where $L_{l,v}$ does not transfer ([Ba1], Lemma 5.2) and the contribution of $l$ on the $G$ side is zero (Prop.\ref{inductionjl} (a) and because $M(s,0)$ is a scalar). If $l$ divides $n$ but $\rho_v$ corresponds to zero (Section \ref{localjl}) for some $v$, then $\tr(\ind(\underbrace{\rho_v\otimes\rho_v\otimes ...\otimes\rho_v}_{l\ {\rm times}})(f_v))=0$ again by Prop.\ref{inductionjl} (a). So the corresponding term is zero (as $M(s,0)$ is a scalar). Using the induction step (b), we see that $\rho_v$ corresponds to zero for some $v$ if and only if $\rho$ is not in the image of ${\bf G}$. 

So we have the equality

$$\sum_{\rho\in DS_{nd}}\tr\rho(f)\ +\ \sum_{l|n,\ l\neq 1}\frac{1}{l^2}\sum_{\rho'\in DS'_{\frac{n}{l}}}\lambda_{{\bf G}(\rho')} 
\tr(\ind({\bf G}(\rho')\otimes {\bf G}(\rho')\otimes ... \otimes {\bf G}(\rho'))(f))$$
$$=\ \sum_{\rho'\in DS'_{n}} m_{\rho'} \tr\rho'(f)\ +\ \sum_{l|n,\ l\neq 1}\frac{1}{l^2}\sum_{\rho'\in DS'_{\frac{n}{l}}} \lambda_{\rho'} \tr(\ind(\rho'\otimes\rho'\otimes ...\otimes\rho')(f'))$$
if $f\lra f'$, where the $\lambda$ are the scalars coming from the intertwining operator $M$ as explained before.

Now, because the local Jacquet-Langlands correspondence commutes with parabolic induction (Proposition \ref{inductionjl} (b)) we have 
$\tr(\ind({\bf G}(\rho')\otimes {\bf G}(\rho')\otimes ... \otimes {\bf G}(\rho'))(f))\ 
=\ \tr(\ind(\rho'\otimes\rho'\otimes ...\otimes\rho')(f'))$ if $f\lra f'$.

If we let $DS_k^{D}$ be the image of ${\bf G}$ and set ${\bf G}(\rho')=\rho$, we get an equality:

$$\sum_{\rho\in DS_{nd}}\tr\rho(f)\ +\ \sum_{l|n,\ l\neq 1}\frac{1}{l^2}\sum_{\rho\in DS^D_{\frac{nd}{l}}}\lambda_{\rho} \tr(\ind(\rho\otimes\rho\otimes ...\otimes\rho)(f))$$
$$=\ \sum_{\rho'\in DS'_{n}} m_{\rho'} \tr\rho'(f')\ +\ \sum_{l|n,\ l\neq 1}\frac{1}{l^2}\sum_{\rho\in DS^D_{\frac{nd}{l}}} 
\lambda_{{\bf G}^{-1}(\rho)} \tr(\ind(\rho\otimes \rho \otimes ...\otimes \rho)(f))$$
if $f\lra f'$.\\

{\bf Let us now prove the case $n=1$.} The previous computations work the same in this case and, as the set of $l\neq 1$ which are divisors of $n=1$ is empty, the equation is:

$$\sum_{\rho\in DS_{nd}}\tr\rho(f)\ =\ \sum_{\rho'\in DS'_{n}} m_{\rho'} \tr\rho'(f),$$
if $f\lra f'$. Now we fix a discrete series $\rho'\in DS'_n$. We let $S$ be the set of places $v$ such that either $v$ is infinite, or  $G'_{1,v}$ is not split, or $\rho'_v$ is not unramified. Then we may separate the representations having the ${}^S$-part isomorphic to ${\rho'}^S$ as explained at Step 4 for the case $n=1$, $d=2$, using the Lemma of the Appendix. We get an equality 
$$\tr(\rho(f))=\sum_{i=1}^k m_{\rho'_i}\tr(\rho'_i(f')),$$
if $f\lra f'$, where $\rho'_i$, $1\leq i\leq k$, are the discrete series in $DS'_1$ which have ${}^S$-part isomorphic to ${\rho'}^S$, and $\rho$ is the only discrete series in $DS_d$ which  have ${}^S$-part isomorphic to ${\rho'}^S$. Such a representation $\rho$ exists because the right hand side is not identically zero by linear independence of characters on $\prod_{v\in S}G'_{1,v}$, and is unique by the strong multiplicity one theorem for $GL_{d}(\a)$. 

Taking $f'_v=f_v=1_{K_v}$ for every $v\notin S$, in order to have $\tr(\rho'_{i,v}(f'_v))=\tr(\rho_v(f_v))=1$ (Chapter 1, Lemma \ref{spherical}), we simplify the formula with all the places not in $S$ and  we get
$$\prod_{v\in S}\tr(\rho_v(f_v))\ =\ \sum_{i=1}^k m_{\rho'_i}\prod_{v\in S}\tr(\rho'_{i,v}(f'_v)),$$
where $f_v=f'_v$ for $v\in S$ such that $G'_1$ splits at $v$ and $f_v\lra f'_v$ for $v\in S$ such that $G'_1$ does not split at $v$.

For every $v\in S$, the representation $\rho_v$ is unitary, and it corresponds to zero or to an irreducible representation $\theta_v$ of $G'_{1,v}$ by the local Jacquet-Langlands correspondence (Theorem \ref{jllocal}). Because on the right hand side we have at least the representation $\rho'_V$, the sum is not identically zero by linear independence of characters on $\prod_{v\in S}G'_{1,v}$, and $\rho_v$ cannot correspond to zero. Then $\otimes_{v\in S}\theta_v$ is an irreducible representation of $\prod_{v\in S}G'_{1,v}$ which satisfies $\prod_{v\in S}\tr(\rho_v(f_v))\ =\ \prod_{v\in S}\tr(\theta_v(f_v))$. By linear independence of characters on $\prod_{v\in S}G'_{1,v}$ (we use the the theory of Harish-Chandra on function characters, since some functions have to have support in the regular semisimple set), we have that on the right there is only one representation, with multiplicity one, and it is isomorphic to $\otimes_{v\in S}\theta_v$. But that representation can be but our $\rho'_V$. Finally $\rho'_v\simeq \theta_v$, so $\rho_v$ corresponds to $\rho'_v$ for all $v\in S$. We already selected $\rho$ in such way that $\rho_v\simeq \rho'_v$ for every $v\notin S$. So we set ${\bf G}(\rho')=\rho$.\\

We prove (c). The previous proof showed that the multiplicity of $\rho'$ was $1$. Let now $\rho''$ be a discrete series of $G'_1$ such that $\rho'_v\simeq \rho''_v$ for all $v\in V$ outside a finite set $T$. Then ${\bf G}(\rho')_v\simeq {\bf G}(\rho'')_v$ for every $v\in V$ such that $G'_1$ splits at $v$ and $v\notin T$. By strong multiplicity one for $G_n$, ${\bf G}(\rho')={\bf G}(\rho'')$. This obviously implies $\rho'_v\simeq \rho''_v$ for every $v\in V$. But the multiplicity of $\rho'$ being one, we have $\rho'=\rho''$. We proved (c).\\

(b) is proved in the following way. Let $\rho$ be a discrete series of $G_n$. If there is a place $v\in V$ such that $\rho_v$ corresponds to zero, it is clear that $\rho$ is not in the image of ${\bf G}$. Let now $\rho$ be such that, for every place $v\in V$, $\rho_v$ corresponds to some representation $\rho'_v$ of $G'_{1,v}$. We start again from the relation
$$\sum_{\rho\in DS_{nd}}\tr\rho(f)\ =\ \sum_{\rho'\in DS'_{n}} m_{\rho'} \tr\rho'(f').$$
Then we fix a finite set of places $S$ containing all the places where $G'_1$ does not split, all the infinite places of $F$ and all the places where $\rho$ is not unramified.
Then we may separate the representations having the ${}^S$-part isomorphic to ${\rho}^S$ as explained at Step 4. We get to an equality 
$$\tr(\rho(f))=\sum_{i=1}^k m_{\rho'_i}\tr(\rho'_i(f')),$$
if $f\lra f'$. Now we use the fact that there exists a couple $f\lra f'$ such that $\tr(\rho(f))\neq 0$. This proves that the equality is not $0=0$ and so there is at least a discrete series $\rho'$ of $G'_1$ on the right. From now the same argument as before shows that ${\bf G}(\rho')=\rho$.\\
\ \\

{\bf Let us prove the theorem for $n$ assuming it is true for $k<n$.} We start with the equation
$$\sum_{\rho\in DS_{nd}}\tr\rho(f)\ +\ \sum_{l|n,\ l\neq 1}
\sum_{\rho\in DS^D_{\frac{nd}{l}}} \frac{\lambda_\rho - \lambda_{{\bf G}^{-1}(\rho)}}{l^2}\ \tr(\ind(\rho\otimes\rho\otimes ...\otimes\rho)(f))$$
$$=\ \sum_{\rho'\in DS'_{n}} m_{\rho'} \tr\rho'(f')$$
if $f\lra f'$. Now we apply the same method as for the proof in the case $n=1$. We fix a discrete series $\rho'$ of $G'_n$, and a finite set of places $S$ 
of $F$ including all the infinite places, the places where $G'_n$ is not split and the places where $\rho'_v$ is not unramified. Then we
separate the representations with ${}^S$-part isomorphic to ${\rho'}^S$. Proposition 4.1 in [Ba1] shows that on the left there is at most one representation with ${}^S$-part isomorphic to ${\rho'}^S$ (the proof uses the unicity of the cuspidal support for representations induced from discrete series of $GL_{nd}(\a)$ 
from [JS] and the classification of the discrete spectrum from [MW1]). So the equality is either of the type of
$$\tr\rho(f)\ =\ \sum_{i=1}^k m_{\rho'_i} \tr\rho'_i(f')$$
and then we continue the proof exactly as in the case $n=1$, or then it is an equality like 
$$\frac{\lambda_{\theta} - \lambda_{{\bf G}^{-1}(\theta)}}{l^2}\ \tr(\ind(\rho\otimes\rho\otimes ...\otimes\rho)(f))\ 
=\ \sum_{i=1}^k m_{\rho'_i} \tr\rho'_i(f'),$$
but we will show that this is impossible.  We remember that $\rho\ =\ {\bf G}(\rho')$ and we get a relation
$$0\ =\ \sum_{i=1}^k m_{\rho'_i} \tr\rho'_i(f') - \beta\  \tr(\ind(\rho'\otimes\rho'\otimes ...\otimes\rho')(f')).$$
where $\beta=\frac{\lambda_{\rho} - \lambda_{{\bf G}^{-1}(\rho)}}{l^2}$.
We simplify with the ${}^S$-part. The trick now is that $l\geq 2$ and $\lambda_{\rho}$ and $\lambda_{{\bf G}^{-1}(\rho)}$ are complex numbers of modulus $1$, so $|\beta|\leq \frac{1}{2}$.
But then, by linear independence of characters on  $G'_{n,S}$, $\beta$ has to be an integer, so it is zero. Then $\sum_{i=1}^k m_{\rho'_i} \tr\rho'_i(f')=0$ which is impossible by the linear independence of the characters on $G'_{n,S}$, since this sum contains at least one term, $\rho'$.

The proof of (b) and (c) is like in the case $n=1$.\\
\ \\

\chapter{Appendix: the spectral simplification}

In this Appendix I recall the proof ([Fl2]) of the classical Lemma which allows a simplification of the equality between the spectral sides of the trace formulae. It is a classical step in the proof of correspondences, and has been used in this paper. At the end of the Appendix, I explain how to adapt  Flath's proof in [FL2] to quasi-split reductive groups (instead of general linear groups).

Let $V$ be a countable set. For almost all $v\in V$ means for all $v\in V$ but a finite set. Let $\{F_v\}_{v\in V}$ be a set of non-Archimedean local fields. For every $v\in V$, let $O_v$ be the ring of integers of $F_v$. 
Let $G_v:=GL_n(F_v)$ and $K_v:=GL_n(O_v)$. On $G_v$ we fix the Haar measure giving volume $1$ to the $K_v$. Let $G$ be the restricted product of the $G_v$, $v\in V$, with respect to $K_v$, $v\in V$.

Let $\H_v^0$ be the set of functions $f:G_v\to \cc$ with compact support and bi-invariant by $K_v$. 
 Let $\H^0$ be the space spanned by functions on $G$ denoted by $f=\otimes_{v\in V}f_v$ with $f_v\in \H_v^0$ for all $v\in V$ and $f_v$ is the characteristic function $1_{K_v}$ of $K_v$ for almost all $v\in V$, and defined by $f((g_v)_{v\in V})\ =\ \prod_{v\in V}f_v(g_v)$. ($\H^0$ is the restricted product of the algebras $\H_v$ with respect to the $1_{K_v}$, where $1_{K_v}$ is the characteristic function of $K_v$, like in Section \ref{restprodalg}.)

Let $J$ be a countable set. For each $j\in J$, for each $v\in V$, let $\pi_{j,v}$ be a unitary spherical smooth irreducible representation of $G_v$.

Assume that, for every couple $j_1,j_2\in J$ such that $j_1\neq j_2$, there is at least one $v\in V$ such that $\pi_{j_1,v}$ and $\pi_{j_2,v}$ are not isomorphic.

\begin{lemme}\label{flath}
Suppose given complex numbers $c_j$, $j\in J$, such that for every $f=\otimes_{v\in V} f_v \in \H^0$, the series
\begin{equation}\label{serie}
\sum_{j\in J}\ c_j\prod_{v\in V}\tr\pi_{j,v}(f_v)
\end{equation}
converges absolutely to zero {\rm (}in $\cc${\rm )}. Then all the $c_j$ are zero.
\end{lemme}
\ \\
{\bf Proof.} 
Let $F$ be a local non-Archimedean field and let $|\ |$ be the normalized norm of $F$.
Every unramified smooth irreducible representation of $GL_n(F)$ is a  sub-quotient of a parabolic induced representation from a character 
$|\ ^.\ |^{z_1}\otimes |\ ^.\ |^{z_2}\otimes ... \otimes |\ ^.\ |^{z_n}$,
($z_i\in \cc$) 
of the diagonal torus with respect to the parabolic subgroup of upper triangular matrices. Conversely, every such induced representation has only
one unramified sub-quotient. See [Car].
Then the set of isomorphism classes of unramified smooth irreducible representations of 
$GL_n(F)$
is in bijection with $\cc^n$ modulo the equivalence relations $R_1$ and $R_2$ defined by 

- $(z_i)_{{1\leq i\leq n}}\cong_{R_1} (z'_i)_{{1\leq i\leq n}}$ if for all $i\in \{1,2,...,n\}$, $z_i- z'_i\in \frac{2\pi}{\ln |u_F|} i{\mathbb Z}$, where $u_F$ is a prime element of $F$, and

- $(z_i)_{{1\leq i\leq n}}\cong_{R_2} (z'_i)_{{1\leq i\leq n}}$ if $(z'_i)_{{1\leq i\leq n}}$ is obtained from  $(z_i)_{{1\leq i\leq n}}$ by a permutation of the components.

Using Tadi\'c classification of unitary representations ([Ta2]) one may show that if such an unramified representation is unitary, then, for all $i$, 
$|Re(z_i)| \leq \frac{n}{2}$, and the set with repetitions $\{Re(z_i)\}_{1\leq i\leq n}$ is symmetrical with respect to zero.
This last condition of symmetry is sufficient for the unramified representation to be hermitian (i.e. isomorphic to the conjugate of its contragredient).

On the set $\mathcal Y$ of n-tuples $(z_1,z_2,...,z_n)\in {\mathbb C}^n$ such that $|Re(z_i)| \leq \frac{n}{2}$ and the set  with repetitions $\{Re(z_i)\}_{1\leq i\leq n}$ is symmetrical with respect to zero put the induced topology by restriction from $\cc^n$. The quotient $\mathcal Y'$ of this space by the equivalence relation $R_1$ is compact. Put the induced topology on $\mathcal Y''$ which is the quotient of $\mathcal Y'$ by the equivalence relation $R_2$.
 The set $\mathcal Y''$ is in natural bijection with the set $\mathcal W$ of isomorphism classes of unramified representations defined by elements in $\mathcal Y$. Consider the induced topology on $\mathcal W$. Then $\mathcal W$ is compact. Every class in $\mathcal W$ is unramified and hermitian; every irreducible unramified unitary class of representations of $GL_n(F)$ is in $\mathcal W$.
 
For each $v\in V$, denote ${\mathcal W}_v$ the analog space of $\mathcal W$ for the group $G_v=GL_n(F_v)$. Let $R:=\prod_{v\in V}{\mathcal W}_v$ be the product space, which is again a compact space. The points of $R$ are collections $\pi=(\pi_v)_{v\in V}$ that have the property that every $\pi_v$ is an (isomorphism class of) irreducible hermitian unramified representation(s) of $G_v$. Not every collection with this property is a point of $R$, but any collection $\pi=(\pi_v)_{v\in V}$ such that every $\pi_v$ is an (isomorphism class of) irreducible unitary unramified representation(s) of $G_v$ is a point of $R$.
For every $f\in \H^0$, define $F_f:R\to \cc$ by the formula $F_f(\pi)=\prod_{v\in V}\tr\pi_v(f_v)$.
The set of all these functions verifies the Stone-Weierstrass conditions:

- it contains the constant functions (consider $f$ equal to a scalar multiple of the characteristic function of $\prod_{v\in V} K_v$); 

- it is stable by complex conjugation because, if $f=\otimes_{v\in V}f_v\in \H^0$, if $f^*(g):=f(g^{-1})$, then $F_{f^*}=\overline{F_f}$ (indeed, 
if $\pi_{v}$ is hermitian, then $\tr\pi_{v}(f_v^*)=\overline{\tr\pi_{v}(f_v)}$); 

- it separates points, because if $\pi_{j_1,v}$ is not isomorphic to $\pi_{j_2,v}$, then there is a function in $\H_v^0$ such that 
$\tr\pi_{j_1,v}(f)\neq  \tr\pi_{j_2,v}(f)$ (particular case of [Be], Corollary 3.9 for example).\\

Hence, the set of all functions $F_f$, $f\in \H^0$,
is a dense subset of the set $C(R)$ of continuous functions on
$R$. Now, putting $f=1$ in formula (\ref{serie}), the absolute 
convergence implies that $\sum_J |c_j|$ converges. Then choose $u$ in $J$
such that $|c_u|$ is maximal. Suppose $|c_u|\neq 0$.
Choose a finite subset $J_0$ of $J$ such that 
$\sum_{J\bc J_0} |c_j|<\frac{|c_u|}{4}$. 
The density result implies  that we may find a function $f$ such that

- $|\mathrm{tr}\pi_j(f)|<2$ for all $j\in J$,

- $|\mathrm{tr}\pi_u(f)|>1$ and

- $|\mathrm{tr}\pi(f)|<\frac{|c_u|}{2|J_0|}$ for all $\pi\in J_0\bc\{u\}$.

Then 
$$|\sum_{j\in J\bc\{u\}}c_j\mathrm{tr}\pi_j(f)|\leq
\sum_{j\in J\bc\{u\}}|c_j||\mathrm{tr}\pi_j(f)|<$$
$$<(|J_0|-1)
\frac{|c_u|}{2|J_0|}
+2\frac{|c_u|}{4}<|c_u|<|c_u\mathrm{tr}\pi_u(f)|$$ 
which contradicts the hypothesis. Hence $|c_u|$ must be $0$.\qed
\ \\
\ \\

Here is the generalization of this technique for other groups. Let $G$ be a connected reductive quasi-split group over a local non archimedean field, let $A$ be a maximal split torus in $G$, $W$ the Weyl group of $A$, $M$ the centralizer of $A$ (so $W$ is the quotient of the normalizer of $A$ by $M$). Then $M$ is a minimal Levi subgroup of $G$. Let $P$ be a parabolic subgroup of $G$ with Levi subgroup $L$. Let $K$ be a maximal compact subgroup of $G$ in good position with respect to $P$. Assume that the Hecke algebra $\H^0$ of left and right $K$-invariant functions over $G$ with compact support is commutative. I will define a compact set ${\mathcal W}$ with the same properties as before, namely ${\mathcal W}$ parametrizes a set of irreducible hermitian unramified representations of $G$ containing all the irreducible unitary unramified representations of $G$. It is known that the set $X$ of unramified characters of $M$ (here unramified means trivial on $M\cap K$) form a complex variety endowed with the obvious action of $W$.

According to Cartier's paper [Car] (I used his notation), for every $\chi\in X$, the induced representation $\ind_P^G\chi$ has a unique irreducible quotient which is an unramified representation of $G$, denoted by $I(\chi)$. Every irreducible unramified representation of $G$ is obtained this way. If $\chi, \chi'\in X$, then $I(\chi)$ and $I(\chi')$ are isomorphic if and only if there is $w\in W$ such that $\chi'=w\chi$. So the set of irreducible unramified representations of $G$ is parametrized by $X/W$. 

Let $X_u$ be the set of characters $\chi\in X$ such that $\ind_P^G\chi$ contains at least one irreducible quotient which has bounded coefficients. Tadi\'c shows, [Ta3] Theorem 7.1, that the set $X_u$ is relatively compact in $X$. It is obvious that $X_u$ is stable by $W$. So the closure $\overline{X_u}$ is compact and stable by $W$.

A representation is hermitian if it is isomorphic to the complex conjugated of its contragredient. The contragredient representation of a character $\chi$ is $\overline{\chi(g^{-1})}$. We also know that the complex conjugation and the contragredient functor commute with parabolic induction. So using the classification, it follows that $I(\chi)$ is hermitian if and only if there is $w\in W$ such that $\overline{\chi(g^{-1})}=w\chi(g)$. If $w\in W$, then the set $Y_w$ of $\chi\in\overline{X_u}$ such that $\overline{\chi(g^{-1})}=w\chi(g)$ is closed, so it is a compact set. $W$ is finite, so $Y:=\cup_{w\in W} Y_w$ is compact. The set ${\mathcal W}:=Y/W$ has the required properties.

\chapter{Bibliography}

[Ar1] J.Arthur, The invariant trace formula, II, Global Theory, {\it J. Am. Math. Soc.} {\bf 1}(3) (1988), 501-554.\\

[Ar2], J.Arthur, An introduction to the trace formula, in {\it Harmonic analysis, the trace formula, and Shimura varieties}, Clay Math. Proc., 4, Providence, R.I.: American Mathematical Society, pp. 1-263.\\

[AC] J.Arthur, L.Clozel, {\it Simple algebras, base change, and the advanced theory of the trace formula}, Annals of Mathematics Studies, 120, Princeton University Press, ISBN 978-0-691-08517-3.\\

[Ba1] A.I.Badulescu, Global Jacquet-Langlands correspondence, multiplicity on and classification of automorphic representations, {\it Inventiones mathematicae} 172(2) (2008), 383-438. With an Appendix by Neven Grbac.\\

[Ba2] A.I.Badulescu,  Jacquet-Langlands et unitarisabilité, J. Inst. Math. Jussieu 6 (2007), no 3, p. 349-379.\\

[Ba3] A.I.Badulescu, Un théorème de finitude. With an appendix by Paul Broussous. Compositio Math. 132 (2002), no. 2, 177-190.\\

[BaRe] A.I.Badulescu, D.Renard,  Unitary dual of $GL_n$ at archimedean places and global Jacquet-Langlands correspondence. Compositio Math 146, vol. 5 (2010), p. 1115-1164.\\

[BaRo] A.I.Badulescu, Ph.Roche, Global Jacquet-Langlands correspondence for division algebra in characteristic $p$, Int.Math.Res.Not. (2017) (7), p. 2172-2206.\\

[Bar] E.M.Baruch, A proof of Kirillov's conjecture, Annals of Mathematics, 158 (2003), 207-252.\\

[Be1] J.N.Bernstein, $P$-invariant distributions on $GL(N)$ and the classification of unitary representations of $GL(N)$ (non-Archimedean case). In {\it Lie Groups and Representations II}, Lect. Notes Math., vol. 1041, Springer (1983).\\

[Be2] J.Bernstein, P.Deligne, Le "centre" de Bernstein,
In {\it Representations des groups redutifs sur un corps local,
Traveaux en cours} (P.Deligne ed.), Hermann, Paris, 1-32 (1984).\\

[Bou1] Bourbaki, Nicolas Integration. I. Chapters 1-6. Translated from the 1959, 1965 and 1967 French originals by Sterling K. Berberian. Elements of Mathematics (Berlin). Springer-Verlag, Berlin, 2004. xvi+472 pp. ISBN: 3-540-41129-1.\\

[Bou2] Bourbaki, Nicolas Integration. II. Chapters 7-9. Translated from the 1963 and 1969 French originals by Sterling K. Berberian. Elements of Mathematics (Berlin). Springer-Verlag, Berlin, 2004. viii+326 pp. ISBN: 3-540-20585-3.\\

[Bou3] Bourbaki Nicolas Algèbre, Chapitre 8. Springer-Verlag 2012.\\

[Bl] D.Blasius, On multiplicities for SL(n), Israel Journal of Mathematics, 88 (1): 237-251.\\

[Bo] A.Borel, Some finiteness properties of adele groups over number fields. Inst. Hautes Études Sci. Publ. Math. No. 16 (1963) 5-30.\\ 

[Booh] J.Booher, The trace formula for compact quotients, with an appendix by B.Conrad, http://math.stanford.edu/~conrad/JLseminar/Notes/L12.pdf\\

[Bu] D.Bump, Automorphic forms and representations. Cambridge Studies in Advanced Mathematics, 55. Cambridge University Press, Cambridge, 1997. xiv+574 pp. ISBN: 0-521-55098-X.\\

[BHLS] A.I.Badulescu, G.Henniart, B.Lemaire et V.Sécherre, Sur le dual unitaire de GL(r,D), Amer. J. Math. 132 (2010), no. 5, 1365-1396.\\

[BJ] A.Borel, H.Jacquet, Automorphic forms and automorphic representations,  in {\it Automorphic forms, representations and $L$-functions, Part 1, Proc. of Symp. in Pure Math. 33}, AMS, 1979, 189-202.\\

[BLM] A.I.Badulescu, E.Lapid, A.Minguez, Une condition suffisante pour l'irr\'eductibilit\'e d'une induite parabolique de $GL(m,D)$,
Annales de l'Institut Fourier, Tome 63 (2013) no. 6, p. 2239-2266.\\

[Car] P. Cartier,   Representations of p-adic groups: a survey, in   {\it Automorphic forms, representations and $L$-functions}, Part 1, Proc. of Symp. in Pure Math. 33, AMS, 1979, 111-155.\\ 

[Cas] W.Casselman, Introduction to the theory of admissible representations of reductive $p$-adic groups, preprint.\\

[Cl] L.Clozel, Théorème d'Atiyah-Bott pour les variétés $p$-adiques et caractères des groupes réductifs,
Mémoires de la Société Mathématique de France, Volume 15, 1984, page 39-64.\\

[De] T.Deng, Induction Parabolique et G\'eométrie des Vari\'et\'es Orbitales pour $GL_n$, PhD thesis, 2016, University Paris 13.\\

[Del] P.Delorme, Infinitesimal Character and Distribution Character of Representations of Reductive Lie Groups, Porceedings of Symposia in Pure Mathematics, Vo. 61 (1997), pp. 73-81.\\ 

[DKV] P.Deligne, D.Kazhdan, M.-F.Vign\'eras, Repr\'esentations des alg\`ebres centrals simples $p$-adiques. In: Repr\'esentations des gropues r\'eductifs sur un corps local, 33-118, Hermann, Paris (1984).\\

[Fl1] D.Flath, Decomposition of representations into tensor products, in {\it Automorphic forms, representations and $L$-functions}, Part 1, Proc. of Symp. in Pure Math. 33, AMS, 1979, 179-184.\\ 

[Fl2] D.Flath, A comparison for the automorphic representations of $GL(3)$ and its twisted forms, Pac.J.Math. 97, 373-402 (1981).\\

[FK] Y.Z.Flicker, D.A.Kazhdan, Metaplectic correspondence. Inst. Hautes Études Sci. Publ. Math. No. 64 (1986), 53-110.\\

[GGP-S] I.M.Gelfand, M.Graev, I.I.Piatetskii-Shapiro, {\it Representation theory and automorphic functions}, W.B. Saunders Company, Philadelphia 1969.\\

[Go1] R.Godement, The spectral decomposition of cusp-forms. in {\it Algebraic Groups and Discontinuous Subgroups} (Proc. Sympos. Pure Math., Boulder,
Colo., 1965) pp. 225-234 Amer. Math. Soc., Providence.\\

[Go2] R.Godement, Domaines fondamentaux des groupes arithmétiques. (French) 1964 Séminaire Bourbaki, 1962/63. Fasc. 3, No. 257 25 pp.\\

[GP-S] I.M.Gelfand, I.I.Pjateckii-Sapiro, Automorphic functions and the theory of representations. (Russian) Trudy Moskov. Mat. Obsc. 12 1963 389-412.\\

[H-C1] Harish-Chandra
{\it Automorphic forms on semisimple Lie groups}, 
Notes by J. G. M. Mars. Lecture Notes in Mathematics, No. 62 Springer-Verlag, Berlin-New York 1968 x+138 pp. \\

[H-C2] Harish-Chandra, A submersion principle and its applications, Proc. Indian Acad. Sci. 90 (1981) 95-102.\\

[He] G.Henniart, La conjecture de Langlands locale pour GL(3). (French) [The local Langlands conjecture for GL(3)] Mém. Soc. Math. France (N.S.) No. 11-12 (1984), 186 pp.\\

[JL] H.Jacquet, R.P.Langlands, {\it Automorphic forms on $GL(2)$} Lecture Notes in Mathematics, Vol. 114. Springer-Verlag, Berlin, 1970.\\

[JS] H.Jacquet, J.A.Shalika, On Euler products and the classification of automorphic forms II, Amer. J. Math. 103 (1981), no. 4, 777-815.\\

[Kn1] A.W.Knapp, Representation theory of semisimple groups. An overview based on examples. Reprint of the 1986 original. Princeton Landmarks in Mathematics. Princeton University Press, Princeton, NJ, 2001. xx+773 pp. ISBN: 0-691-09089-0.\\

[Kn2] A.W.Knapp, Advanced algebra, Cornerstones, Birkhäuser, Boston, 2007, ISBN: 0-8176-4522-5,
730+xxiv pages.\\

[KnVo] A.Knapp, D.A.Vogan, {\it Cohomological Induction and Unitary Representations}, Princeton University Press, Princeton (1995).\\

[Lau] G.Laumon, {\it Cohomology of Drinfeld Modular Varieties}, Part I, Part II. Cambridge University Press, 1996.\\


[MW1] C.Moeglin, J.-L.Waldspurger, Le spectre résiduel de $GL_n$, {\it Ann.Sci.Ecole.Norm.Sup.} (4) 22 (1989), 605-674.\\





[P-S] I.Piatetski-Shapiro, Multiplicity one theorems, in {\it Automorphic forms, representations and $L$-functions}, Part 1, Proc. of Symp. in Pure Math. 33, AMS, 1979, 209-212.\\ 

[PR] V.Platonov, A.Rapinchuk, Algebraic groups and number theory. Translated from the 1991 Russian original by Rachel Rowen. Pure and Applied Mathematics, 139. Academic Press, Inc., Boston, MA, 1994. xii+614 pp. ISBN: 0-12-558180-7.\\

[Re] D.Renard, Repr\'esentations des groupes r\'eductifs $p$-adiques, Cours sp\'ecialis\'es, 17, SMF.\\

[Ro1] J.Rogawski, Representations of ${\rm GL}(n)$ and division algebras over a $p$-adic field, Duke Math. J. 50 (1983), no. 1, 161-196.\\

[Ro2] J.Rogawski, {\it Automorphic Representation of Unitary Groups in Three Variables}, Annals of Mathematical Studies 123, Princeton University Press, Princeton 1990.\\

[RV] D. Ramakrishnan, R. Valenza, {\it Fourier analysis on number fields}. Graduate Texts in Mathematics, 186. Springer-Verlag, New York, 1999.\\

[Se] V.Sécherre, Proof of the Tadi\'c Conjecture U0 on the unitary dual of GL(m,D), J. Reine Angew. Math. 626 (2009), 187-204.\\

[Sh] J.A.Shalika, The multiplicity one theorem for $GL_n$, Annals of Mathematics. Second Series, 100: 171-193.\\

[Ta1] M.Tadi\'c, Representation theory of GL(n) over a p-adic division algebra and unitarity in
the Jacquet-Langlands correspondence, Pacific J. Math. 223 (2006), no. 1, p. 167-200.\\

[Ta2] M.Tadi\'c, Classification of unitary representations in irreducible representations of general
linear group (non-Archimedean case), Ann.Sci.\'Ecole Norm.Sup. (4) 19 (1986), no. 3, p. 335-382.\\

[Ta3] M.Tadi\'c, Geometry of dual spaces of reductive groups (non archimedean case), Journal d'Analyse Math\'ematique, vol. 51, 1988, pages 139-181.\\

[Ti] J.Tits, Reductive groups over local fields, in {\it Automorphic forms, representations and $L$-functions, Part 1, Proc. of Symp. in Pure Math. 33}, AMS, 1979, 29-71.\\

[Vo] D.A.Vogan, The unitary dual of $GL(n)$ over an Archimedean field, Invent.Math. 83 (1986), 449-505.\\

[We] A.Weil, {\it Basic number theory}, Reprint of the second (1973) edition. Classics
in  Mathematics.  Springer-Verlag,  Berlin,  1995.  xviii+315.\\

[DL] M.Duflo, J.-P.Labesse Sur la formule des traces de Selberg, Ann.Sci.ENS 4, 2, (1971), pp. 193-284.\\

[CL2] L.Clozel, Invariant Harmonic Analysis on the Schwartz Space of a Reductive p-adic Group, in {\it Harmonic Analysis on Reductive Groups}, Progress in Mathematics (101), Birkhäuser, Boston, (1991), pp 101-121.\\

[GJ] S.Gelbart, H.Jacquet, Forms on $GL(2)$ from the automorphic point of view, pp.213-252.\\

[La] R.P.Langlands, {\it On the Functional Equations Satisfied by Eisenstein Series}, Lecture Notes in Math., Vol. 544, Springer-Verlag, Berlin-Heidelberg-New York, 1976.\\

[LL] J.-P.Labesse, R.P.Lannglands, L-indistinguishability for SL(2), Canadian J. Math. 31 (1979) pp.726-785. 

[Serre] J.-P. Serre, Corps locaux, 3 edition, Herman, 1080.

\end{document}